\documentclass[11pt]{article}

\usepackage{amsfonts,amsmath}

\usepackage[dvips]{graphicx}
\usepackage{color}
\usepackage{float}
\usepackage{epsfig} 
\usepackage{fullpage}
\newtheorem{theorem}{Theorem}[section]
\newtheorem{cor}[theorem]{\bf{Corollary}}
\newtheorem{lem}[theorem]{\bf{Lemma}}

\newtheorem{prop}[theorem]{\bf{Proposition}}
\newtheorem{df}[theorem]{\bf{Definition}}
\newtheorem{ex}[theorem]{\bf{Example}} 

\newtheorem{rem}[theorem]{\bf{Remark}}

\usepackage{epstopdf}% To incorporate .eps illustrations using PDFLaTeX, etc.
\usepackage{subfigure}% Support for small, `sub' figures and tables

\title{Joint dynamic probabilistic constraints with \\ projected linear decision rules}
\author{Vincent Guigues\\ FGV/EMAp,\\ 22250-900 Rio de Janeiro, Brazil\\ {\tt vguigues@fgv.br} \and Ren\'e Henrion\\Weierstrass Institute Berlin\\10117 Berlin, Germany\\henrion@wias-berlin.de}
\date{}

\begin{document}
\maketitle

\begin{abstract}
We consider  multistage  stochastic  linear  optimization  problems 
combining joint dynamic probabilistic constraints with
hard constraints. We develop a method for projecting decision rules 
onto hard constraints of wait-and-see type. We establish the relation between the original (infinite-dimensional) problem and approximating problems working with projections 
from different subclasses of decision policies.
Considering the subclass of linear decision rules and a generalized linear model for the underlying stochastic process
with noises that are Gaussian or truncated Gaussian, we show that the value and gradient of the objective and constraint functions
of the approximating problems can be computed analytically.
\end{abstract}



\section{Introduction}
%{\color{red} \`a compl\'eter}

\noindent 

Probabilistic constraints were introduced some fifty
years ago under the name 'chance constraints' by Charnes and
Cooper \cite{charnes}.
A probabilistic constraint is an inequality 
\begin{equation}
\mathbb{P}\left( g(x,\xi )\leq 0\right) \geq p,  \label{chancecons}
\end{equation}%
where $g$ is a mapping defining a random inequality system, $x$ is a decision vector, and $\xi $ is a random vector living on a probability space $(\Omega,\mathcal{A},\mathbb{P})$.
The meaning of (\ref{chancecons}) is the following:
a decision $x$ is feasible if and only if the
random inequality system $g(x,\xi )\leq 0$ is satisfied at least with
probability $p\in (0,1]$. Choosing $p$ close to one reflects the wish for robust decisions which can be interpreted in a probabilistic way. 

In the beginning, efforts focussed on finding
explicit deterministic equivalents for (\ref{chancecons}), i.e., on finding
analytical functions such that (\ref{chancecons}) is equivalent with the
inequality $\varphi (x)\geq p$, see \cite{mill} for instance. 
Even if such instances are rare and usually related with special
assumptions, e.g., one-dimensional random variables, individual
probabilistic constraints, or assuming
independent components of the random vector, it has been successfully applied more recently
using Boolean Programming to attack
joint probabilistic constraints with dependent random right-hand sides \cite{LEJE12a, LEJE12b} and later extended 
to stochastic programming problems with joint probabilistic constraints and multi-row random technology matrix \cite{KOLE12}.

A new era in the theoretical and algorithmical treatment of probabilistic
constraints began with the pioneering work by Pr\'{e}kopa in the early
seventies, when the theory of log-concave probabilty measures allowed to
derive the convexity of feasible decisions induced by a large class of
probabilistic constraints (\ref{chancecons}). Along with bounding and
simulation techniques outperforming crude Monte Carlo approaches, this paved
the way for applying efficient methods from convex optimization for the
numerical solution of probabilistically constrained optimization problems.
The monograph \cite{prek} is still a standard reference in
this area. 

Another breakthrough in this direction happened in the early
nineties and was related with efficient codes for numerical
integration of multivariate normal and t-probabilities due to
Genz \cite{genz1}. These codes are to the best of our knowledge
the best performing ones in this area, up until now. For a recent
survey on this topic, we refer to the monograph \cite{genz}. Along
with a reduction technique which allows us to lead back analytically
the computation of gradients to the computation of values in
(\ref{chancecons}), these codes may be used for solving
probabilistically constrained optimization problems in meaningful
dimension of up to a few hundred (as far as the random vector is
concerned). For some recent applications in energy management, we
refer to \cite{ackooij4} and \cite{ackooij2}.

Alternative solution methods rely on convex approximations of the chance constraints,
see for instance \cite{nemsha06} (where Bernstein approximations are used) and \cite{chenglisser},
and on the scenario approach to build computationally tractable approximations as in \cite{cala1}, \cite{cala2}, \cite{farias}.

Applications of probabilistic constraints are abundant in engineering and
finance (for an overview on the
theory, numerics and applications of probabilistic constraints, we refer to,
e.g., \cite{shap-rusc-book}, \cite{prek}, and \cite{prek4}). Within
engineering, power management problems are dominating as far as
probabilistic constraints are concerned. In particular, hydro reservoir management is a
fruitful instance for this class of optimization problems. We
may refer to the basic monograph \cite{loucks} and to some exemplary work in
this field (\cite{chat}, \cite{duran}, \cite{edi}, \cite{gs2010}, \cite%
{loia}, \cite{mor}, \cite{prek0}, \cite{prek1}).

In many applications, the decision $x$ has to be taken before the realization of the 
random parameter $\xi$ is observed ('here-and-now decisions').
However, decisions often depend on time, i.e., the vector $x$ represents a discrete decision process. In such case, the 'here-and-now' setting of 
(\ref{chancecons}) means that decisions for the whole time period are taken prior to observing the random parameter, which is now a discrete stochastic process. Then
inequality \eqref{chancecons} represents a static probabilistic constraint because the decision process does not take into account the gain of information over time while observing the random process. To overcome 
this deficiency, one may pass from a decision vector $x=(x_1,\ldots ,x_T)$ to a closed-loop decision policy 
\begin{equation}\label{policyform}
x=(x_1,x_2(\xi_1),x_3(\xi_1,\xi_2),\ldots , x_T(\xi_1,\ldots\xi_{T-1}))
\end{equation}
each component of which represents a function of previously observed values of the random process for a given time. 
A simple way to compute a closed-loop strategy
is the application of a rolling horizon policy which at any time of the
horizon hedges against future uncertainty conditional to past realizations
of the random process (see, e.g., \cite{prek0}, \cite{prek1}, \cite{gs2010}, \cite{gsorl2012}, \cite{gs2012}). Only the obtained optimal decision for the next time
step is applied in reality. Another possibility consists in computing the policy at the beginning of the optimization period
plugging \eqref{policyform} into \eqref{chancecons} and \eqref{chancecons}
becomes a dynamic probabilistic constraint now acting on a variable $x$ from an infinite-dimensional space. 

In this setting, in order to return to a numerically tractable problem in finite dimensions, the decision policies are often parameterized, the most common approach being the introduction of linear decision rules, 
i.e.,
$x_i(\xi)=A\xi +b$ for appropriate $A,b$ which now become the finite-dimensional substitutes for the originally infinite-dimensional variables. This strategy has been introduced to 
probabilistically constrained hydro reservoir problems as early as 1969  \cite{revelle}. It was used there (and in subsequent publications) in the context of so-called 
individual probabilistic constraints where each component of the given random inequality system is individually turned into a probabilistic constraint:
\[
\mathbb{P}\left( g_i(x,\xi )\leq 0\right) \geq p\quad (i=1,\ldots ,m).
\]
The big advantage of such individual constraints is that - in case the component $g_i(x,\xi)$ is separable with respect to $\xi$ - they are easily converted 
into explicit constraints via quantiles. In particular, if $g$ happens to be a linear mapping and the objective is linear too, then all one has to do to solve such a probabilistic optimization problem is to apply linear programming. It is well known, however, that the 
probability level $p$ chosen in an 
individual model may by far not correspond to the level in a joint model, given by  \eqref{chancecons}, where the probability is taken over the entire inequality system. In \cite{ackooij4} a hydro reservoir 
problem is presented where at an optimal release policy the level constraints are satisfied in each time interval with probability 90\% individually, whereas the probability of keeping the level constraints through the whole time period is as low as 32\%. This observation strongly suggests to deal with the joint model \eqref{chancecons} 
albeit much more difficult to treat algorithmically.

Joint probabilistic constraints in the closed-loop sense discussed above have been investigated in \cite{henrion} again 
in the context of a reservoir problem. Here a highly flexible piecewise constant approximation of decision policies $x(\xi)$ was considered and it turned out that the optimal policies of 
the given problem were definitely not linear. However, a sufficiently fine piecewise approximation requires a big computational effort and limits the applicability of the model to a few 
time stages like three or four. Therefore, picking up again the idea of parameterized (in particular, linear) decision rules but now in the context of joint constraints appears to be reasonable.

Other authors embed optimization problems with dynamic probabilistic constraints into a dynamic programming scheme of optimal control, however, typically 
imposing simplifications with regard to the joint system of constraints like the assumption of independent components, or of a discrete 
distribution (scenarios) or of an individualized (via Boole-Bonferroni inequality) surrogate model (e.g., \cite{carpentier,ono}).

The aim of the current paper is to discuss several modeling issues in the context of dynamic probabilistic constraints putting the emphasis on
\begin{itemize}
\item joint probabilistic constraints as in \eqref{chancecons};
\item continuous multivariate distributions of the random vector (in particular, Gaussian) with typically correlated components;
\item parameterized decision rules (in particular, linear and projected linear ones); and
\item mixed probabilistic and hard (almost sure) constraints.

\end{itemize}
We do not intend to investigate the so-called time consistent models for dynamic probabilistic constraints as it was done, for instance, in \cite{carpentier}. 
This issue has been considered so far in the framework of Dynamic Programming, where the assumption of the random vector having independent components 
is paramount, e.g., \cite{alaiscarplara}. Moreover, typically, a discrete distribution is assumed for numerical analysis. As pointed out above, we are interested here 
in continuously distributed distributions with potentially correlated variables. Though it seems possible to establish time consistent models for dynamic chance 
constraints under multivariate Gaussian distribution, this issue would complicate the analysis we have in mind here and is yet to be explored in future research.

Moreover, the focus of this paper is not 
to develop a new algorithm neither the study of a concrete application, although a simple hydro reservoir problem will guide us as an illustration. Our idea is rather to provide a modeling framework taking 
into account the items listed above and yielding a link to algorithmic approaches for static probabilistic constraints. 
The latter have been successfully dealt with numerically in the context of linear probabilistic constraints under multivariate Gaussian (and Gaussian-like) distribution 
(see, e.g., \cite{prek,prek0,ackooij4,ackooij2}).

The paper is organized as follows: Section 2 presents a general linear multistage problem with probabilistic and hard constraints. It describes a method for projecting decision rules 
onto hard constraints of wait-and-see type. It finally establishes the relation between the original (infinite-dimensional) problem and approximating problems working with projections 
from different subclasses of decision policies. These subclasses are kept very general in this section while they are specialized to linear decision rules in Section 3. In that 
same section the probabilistic time series model we intend to use for the discrete stochastic process is made precise. It is clarified, how the objective, the probabilistic constraint and the hard constraints 
look like under this probabilistic model and the assumed linear decision rules. Finally, Section 4 explicitly develops the shape of general optimization problems introduced in Section 2 when assuming 
multivariate Gaussian and truncated Gaussian models for the discrete process. Advantages and difficulties for the different problems are discussed.
\section{A linear multistage problem with probabilistic constraints}
\subsection{The general model}
\noindent For given $T\in \mathbb{N}$ with $T \geq 2$, we consider a $T$-stage stochastic 
linear
minimization problem with  the following random constraints:
\begin{eqnarray} \label{randopt0}
%&\mbox{minimize}\quad\sum_{t=1}^{T}\langle h_{t},y_{t} \rangle\quad\mbox{subject to}&\nonumber\\
&\sum\limits_{\tau =1}^{t}A_{t, \tau }y_{\tau }+\sum\limits_{\tau
=1}^{t}B_{t, \tau }\xi _{\tau }\leq b_{t},\quad t=1,\ldots ,T.&
\end{eqnarray}
Here, for $t=1,\ldots ,T$, $y_{t}$ are 
$n_t$-dimensional decision vectors, $\xi _{t}$ are 
$M_t$-dimensional random vectors, $A_{t, \tau }$ and $B_{t, \tau }$ are given matrices of orders $(l_t,n_\tau)$ and $(l_t,M_\tau)$, respectively, and
$b_{t}\in\mathbb{R}^{l_t}$ are given vectors. In what
follows, the index '$t$' will be interpreted as time and $y_{t}$ and $%
\xi _{t}$ represent discrete decision and stochastic processes,
respectively, having finite horizon. In this time-dependent setting, we shall assume that all components of the random process have the same dimension $M_1=\cdots =M_T=:M$.
The joint random vector $\xi =(\xi_1\ldots ,\xi_T)\in\mathbb{R}^{MT}$ is supposed to live in a probability space $\left( \Omega ,\mathcal{A},\mathbb{P}\right) $.
Similarly to traditional multistage stochastic programming, we shall assume
that the decision $y_{t}$ is taken in the beginning of time interval $%
[t,t+1) $ but the random vector $\xi _{t}$ is observed only at the end
of that same interval. Therefore, the realization of $\xi _{t}$ is unknown
at the time one has to decide on $y_{t}$. On the other hand, in order to take into
account the gain of information due to past observations of randomness, the
decision $y_{t}$ is allowed to depend on $\xi _{1:t-1}:=\left( \xi
_{1},\ldots ,\xi _{t-1}\right) $ such that $y_{t}$ is Borel measurable. 
In the following, we will refer to the $y_{t}\left( \xi _{1:t-1}\right), t=1,\ldots,T$,
(including the deterministic first stage decision $y_{1}\left( \xi
_{1:0}\right) :=y_{1}$) as decision policies rather than decision vectors in
order to emphasize their functional character.
Summarizing, we are dealing with the following problem:
\begin{eqnarray} \label{randopt0a}
&\mbox{minimize}\quad \mathbb{E}\sum_{t=1}^{T}\langle h_{t},y_{t}\left( \xi _{1:t-1}\right)\rangle\quad\mbox{subject to}&\nonumber\\
&\sum\limits_{\tau =1}^{t}A_{t, \tau }y_{\tau }\left( \xi _{1:\tau-1}\right)+\sum\limits_{\tau
=1}^{t}B_{t, \tau }\xi _{\tau }\leq b_{t},\quad t=1,\ldots ,T,&
\end{eqnarray}
where $\mathbb{E}$ is the expectation operator and $h_t$ is a deterministic cost vector for stage $t$.
\begin{ex}\label{hydroex1}
As an illustration, we
consider a two-stage problem for the optimal release $y$ of a
hydro-reservoir under stochastic inflow $\xi $. The released water is used
to produce and sell hydro-energy at a price $p$ which is assumed to be known
in advance. Given the two stages, these quantities have components $\xi
=(\xi _{1},\xi _{2})$, $p=(p_{1},p_{2})$, $y=(y_{1},y_{2}(\xi _{1}))$. The
reservoir level is required to stay at both stages between given lower and
upper limits $\ell^{lo}$, $\ell^{up}$, respectively. Finally, the release is
supposed to be bounded by fixed operational limits $y^{lo}$, $y^{up}$,
respectively, for turbining water at both time stages. Denoting by $\ell_{0}$
the initial water level in the reservoir, the 
random cost is given by $-(p_{1}y_{1}+p_{2}y_{2}(\xi _{1}))$
while the random constraints can be written
\begin{equation} \label{hydro}
\begin{array}{c}
\ell^{lo}\leq \ell_{0}+\xi _{1}-y_{1}\leq \ell^{up} \\
\ell^{lo}\leq \ell_{0}+\xi _{1}+\xi _{2}-y_{1}-y_{2}(\xi _{1})\leq \ell^{up} \\
y^{lo}\leq y_{1}\leq y^{up} \\
y^{lo}\leq y_{2}(\xi _{1})\leq y^{up}.
\end{array}
\end{equation}
It is easy to see that this is a special instance of problem (\ref{randopt0a})
with data
\begin{eqnarray*}
&h:=-p,\,A_{1, 1}:=A_{2, 2}:=\left( 
\begin{array}{r}
-1 \\ 
1 \\ 
1 \\ 
-1
\end{array}
\right) ,\,A_{2, 1}:=\left( 
\begin{array}{r}
-1 \\ 
1 \\ 
0 \\ 
0
\end{array}
\right), & \\
&B_{1, 1}:=B_{2, 1}:=B_{2, 2}:=\left( 
\begin{array}{r}
1 \\ 
-1 \\ 
0 \\ 
0
\end{array}
\right) ,\,b_{1}:=b_{2}:=\left( 
\begin{array}{c}
\ell^{up}-\ell_{0} \\ 
\ell_{0}-\ell^{lo} \\ 
y^{up} \\ 
-y^{lo}
\end{array}
\right) .&
\end{eqnarray*}
\end{ex}

\bigskip\noindent
%In order to make minimization in the random optimization problem \eqref{randopt0a} a meaningful concept, the random objective is
%usually replaced by its expected value or some risk measure operating on it.
%We shall consider the risk-neutral case using expected values here. 
As far as the
constraints are concerned, satisfying them in expectation only, would result
in decisions leading to frequent violation of constraints which is not
desirable for a stable operation say of technological equipment, etc. At the
other extreme, constraints could be required to hold almost surely, thus
yielding very robust decisions avoiding violation of constraints with
probability one. In that case, %(\ref{randopt})
we obtain the
well-defined optimization problem%
\begin{eqnarray}\label{detopt}
&\mbox{minimize }\mathbb{E}\sum_{t=1}^{T}\langle h_{t},y_{t}\left( \xi
_{1:t-1}\right) \rangle \quad \mbox{subject to}&  \notag \\
&\sum\limits_{\tau =1}^{t}A_{t, \tau }y_{\tau }\left( \xi _{1:\tau -1}\right)
+\sum\limits_{\tau =1}^{t}B_{t, \tau }\xi _{\tau }\leq b_{t}\quad t=1,\ldots
,T,&\quad \mathbb{P}\mbox{-almost surely}.  
\end{eqnarray}%
If in the constraints of (\ref{detopt}) one had that $B_{T, T}=0$, then the
last component $\xi _{T}$ of the random process would not enter the
constraints and (\ref{detopt}) would represent a conventional multistage
stochastic linear program. Note, however, that $B_{2, 2}\neq 0$ in the
two-stage problem (\ref{hydro}) and so the random inflow $\xi _{2}$ observed
only after taking the last decision $y_{2}(\xi _{1})$ plays a role in some
of the (level) constraints. In such cases, insisting on almost sure
satisfaction of constraints may be impossible in particular for unbounded
random distributions. In (\ref{hydro}), for instance, no matter what has
been observed ($\xi _{1}$) or decided on ($y_{1}$,$y_{2}(\xi _{1})$) until
the beginning of the second time interval, the last unknown inflow $\xi _{2}$
could always be large enough to eventually violate the upper-level
constraint 
\begin{equation*}
\ell_{0}+\xi _{1}+\xi _{2}-y_{1}-y_{2}(\xi _{1})\leq \ell^{up}.
\end{equation*}%
Therefore, one has to look for alternative models for such constraints
leaving the possibility of a 'controlled' violation. These observations lead
us to distinguish in (\ref{detopt}) between \textit{hard constraints} which
have to be satisfied almost surely for physical or logical reasons and 
\textit{soft constraints} which can be dealt with in a more flexible way. A
typical example for hard constraints are the lower and upper limits for the
amounts of turbined water ($y^{lo}\leq y_{1},y_{2}(\xi _{1})\leq y^{up}$) in
(\ref{hydro}): there is no turbining beyond the given operational limits
just for physical reasons.

On the other hand, the reservoir level constraints could be considered to be
soft ones. Suppose, for instance, that $\ell^{lo}$ in (\ref{hydro}) represents
the physical lower limit of the reservoir below which no water is released
and turbined. Then, a violation of the lower-level constraint can never
happen and so the corresponding two inequalities can be removed from (\ref%
{hydro}). Doing so, one has to take into account, however, that not the
total amount of the release policies $y_{1}$ and $y_{2}(\xi _{1})$,
respectively, can be turbined and sold at the given prices but only the part
not violating the lower-level constraint, i.e., $\min \{y_{1}, \ell_{0}+\xi
_{1}-\ell^{lo}\}$ in the first stage and $\min \{y_{2}(\xi _{1}), \ell_{0}+\xi
_{1}+\xi _{2}-l^{lo}-y_{1}\}$ in the second stage. This means that the
original profits $p_{1}y_{1}$ and $p_{2}y_{2}(\xi _{1})$ at the two stages
have to be reduced by the amounts $p_{1}\left( y_{1}-\ell_{0}-\xi
_{1}+\ell^{lo}\right) _{+}$ and $p_{2}\left( y_{2}(\xi _{1})-\ell_{0}-\xi _{1}-\xi
_{2}+\ell^{lo}+y_{1}\right) _{+}$, respectively, where the lower index '+' as usual represents the component-wise maximum of the given expression and zero. In this way, the original
lower-level constraints in (\ref{hydro}) have been removed and compensated
for by appropriate penalty terms in the objective.

Next, suppose that $\ell^{up}$ in (\ref{hydro}) represents some upper limit of
the reservoir which is considerably lower than the physical one and serves
the purpose of keeping a flood reserve. Then we may neither be able to
satisfy this upper limit almost surely (see above) nor to remove it in
exchange for an appropriate penalty. In such cases it is reasonable to
impose a probabilistic constraint instead:%
\begin{equation*}
\mathbb{P}\left( \ell_{0}+\xi _{1}-y_{1}\leq \ell^{up},\,\ell_{0}+\xi _{1}+\xi
_{2}-y_{1}-y_{2}(\xi _{1})\leq \ell^{up}\right) \geq p,
\end{equation*}%
where $p\in (0,1)$ is a specified probability level. Hence, the release
policies $y_{1},y_{2}(\xi _{1})$ are defined to be feasible if the indicated
set of random inequalities is satisfied at least with probability $p$.
Observe that $p=1$ would yield the almost sure constraints again, hence
choosing $p$ close to but smaller than one, offers us the possibility of finding a
feasible release policy while keeping the soft upper-level constraint
in a very robust sense.
\begin{ex}\label{hydroex2}
Taking into account all three kinds of hard and soft constraints in the
(random) hydro reservoir model (\ref{hydro}), one ends up with the following
well-defined optimization problem:
\begin{eqnarray}
&
\begin{array}{l}
\mbox{minimize } \\ 
-\mathbb{E}(p_{1}y_{1}+p_{2}y_{2}(\xi _{1})) \\ 
+\mathbb{E}(p_{1}(y_{1}-\ell_{0}-\xi _{1}+\ell^{lo})_{+}+p_{2}(y_{2}(\xi
_{1})-\ell_{0}-\xi _{1}-\xi _{2}+y_{1}+\ell^{lo})_{+})
\end{array}
&  \label{hydrototal} \\
&&  \notag \\
&\mbox{subject to}&  \notag \\
&&  \notag \\
&\mathbb{P}\left( 
\begin{array}{r}
\ell_{0}+\xi _{1}-y_{1}\leq l^{up} \\ 
\ell_{0}+\xi _{1}+\xi _{2}-y_{1}-y_{2}(\xi _{1})\leq \ell^{up}
\end{array}
\right) \geq p&  \notag \\
&&  \notag \\
&\left. 
\begin{array}{c}
y^{lo}\leq y_{1}\leq y^{up} \\ 
y^{lo}\leq y_{2}(\xi _{1})\leq y^{up}
\end{array}
\right\} \quad \mathbb{P}\mbox{-almost surely.}&  \notag
\end{eqnarray}
Here, the group of soft lower-level constraints has disappeared and entered
the objective as a second penalization term, the group of soft upper-level
constraints (for which no penalization costs are available) has turned into
a probabilistic constraint and the group of hard box constraints is
formulated in the almost sure sense.
\end{ex}

\bigskip\noindent
Applying this strategy to the general random constraints (\ref%
{detopt}), we are led to partition the data matrices and vectors for
\thinspace $t=1,\ldots ,T$, \thinspace\ and \thinspace $\tau =1,\ldots ,t$,%
\thinspace\ as 
\begin{equation*}
A_{t, \tau }=\left( A_{t, \tau }^{(1)},A_{t, \tau }^{(2)},A_{t, \tau }^{(3)}\right)
,\,\,B_{t, \tau }=\left( B_{t, \tau }^{(1)},B_{t, \tau }^{(2)},B_{t, \tau
}^{(3)}\right) ,\,\,b_{t}=\left( b_{t}^{(1)},b_{t}^{(2)},b_{t}^{(3)}\right)
\end{equation*}%
according to penalized soft constraints (upper index (1)), probabilistic
soft constraints (upper index (2)) and almost sure hard constraints (upper
index (3)). Accordingly, (\ref{detopt}) turns into the well-defined
optimization problem 
\begin{eqnarray}
&\mbox{minimize }&  \label{finopt} \\
&\sum\limits_{t=1}^{T}\mathbb{E}\left\{ \langle h_{t},y_{t}\left( \xi
_{1:t-1}\right) \rangle + \left< \mathcal{P}_{t}, \left( \sum\limits_{\tau
=1}^{t}A_{t, \tau }^{(1)} y_{\tau }(\xi _{1:\tau -1})+\sum\limits_{\tau
=1}^{t}B_{t, \tau }^{(1)}\xi _{\tau }-b_{t}^{(1)}\right) _{+} \right> \right\} &  \notag \\
&\mbox{subject to}&  \notag \\
&\mathbb{P}\left( \sum\limits_{\tau =1}^{t}A_{t, \tau }^{(2)}y_{\tau }\left(
\xi _{1:\tau -1}\right) +\sum\limits_{\tau =1}^{t}B_{t, \tau }^{(2)}\xi _{\tau
}\leq b_{t}^{(2)},\quad t=1,\ldots ,T\right) \geq p&  \notag \\
&\sum\limits_{\tau =1}^{t}A_{t, \tau }^{(3)}y_{\tau }\left( \xi _{1:\tau
-1}\right) +\sum\limits_{\tau =1}^{t}B_{t, \tau }^{(3)}\xi _{\tau }\leq
b_{t}^{(3)},\quad t=1,\ldots ,T,\quad \mathbb{P}\mbox{-almost surely.}&  \notag
\end{eqnarray}%
Here, the $\mathcal{P}_{t}\geq 0$ refer to cost vectors penalizing the violation of soft constraints with upper index (1).

\subsection{Projection onto hard constraints of wait-and-see type}\label{hardwait}

\noindent We will refer in (\ref{detopt}) to \textit{%
wait-and-see constraints }if $B_{t, t}=0$ for all $t=1,\ldots ,T$, and to 
\textit{here-and-now constraints} otherwise. The distinction is made
according to whether in the constraint of any stage $t$ there is unobserved
randomness $\xi _{t}$ left or not. For example, in (\ref{hydro}), the first
two inequalities (level constraints) are here-and-now whereas the last two
(operational limits) are wait-and-see. As mentioned earlier, the almost sure
constraints in (\ref{finopt}) do not  have a good chance to be ever satisfied
if $B_{T, T}^{(3)}\neq 0$ and the support of the random distribution is
unbounded. We will  get back to such here-and-now constraints for bounded
support of the random distribution in Section \ref{truncgaussopt}. First, let us deal with
the case where all hard constraints are of wait-and-see type as in (\ref{hydrototal}%
). In this case, owing to $B_{t, t}^{(3)}=0$ for all $t=1,\ldots ,T$, the constraint set of (\ref{finopt}) can be written as%
\begin{align}
& M_{1}:=\Big\{\left( y_{t}\left( \xi _{1:t-1}\right) \right) _{t=1,\ldots ,T}|
& &  \label{m1def} \\
& \mathbb{P}\left( \sum\limits_{\tau =1}^{t}A_{t, \tau }^{(2)}y_{\tau }\left(
\xi _{1:\tau -1}\right) +\sum\limits_{\tau =1}^{t}B_{t, \tau }^{(2)}\xi _{\tau
}\leq b_{t}^{(2)},\quad t=1,\ldots ,T\right) \geq p & &  \notag \\
& \sum\limits_{\tau =1}^{t}A_{t, \tau }^{(3)}y_{\tau }\left( \xi _{1:\tau
-1}\right) +\sum\limits_{\tau =1}^{t-1}B_{t, \tau }^{(3)}\xi _{\tau }\leq
b_{t}^{(3)},\quad t=1,\ldots ,T,\quad \mathbb{P}\mbox{-almost surely} \Big\}. & &  \notag
\end{align}
In the context of numerical solution approaches, one
will usually not work in the infinite-dimensional setting of all Borel
measurable policies but rather with a finite-dimensional approximation which
may be defined by some proper subset $\mathcal{K}$ of policies. Later in
this paper we will deal with the class of \textit{linear decision rules} (see Section \ref{lindec}). The feasible set of (\ref{finopt}) will then become the intersection $M_1\cap\mathcal{K}$ rather than just $M_1$. This intersection may turn out to be very small or even empty thus leading to a poor approximation of the 
infinite-dimensional problem (\ref{finopt}). If, for instance, one of the hard constraints is given as 
$y_2(\xi_1)\in [1,2]$ ($\mathbb{P}$-almost surely) and if, moreover, the class of policies 
is
\[
\mathcal{K}:=\{(y_1,y_2(\xi_1)) |\exists a\in\mathbb{R}:y_2(\xi_1)=a\xi_1\},
\]
then, clearly,  $M_1\cap\mathcal{K}=\emptyset$. One possibility to avoid this kind of 
problem is to operate with projections of policies onto the feasible domain of hard constraints. 

Given a closed convex subset $X$ of a finite-dimensional space, we
denote the uniquely defined projection onto this set by $\pi _{X}$. For $t=1,\ldots ,T$, we introduce the multifunctions
\begin{equation}
X_{t}\left( z_{1:t-1},\xi _{1:t-1}\right) :=\left\{ y|A_{t, t}^{(3)}y(\xi_{1:t-1})\leq
b_{t}^{(3)}-\sum\limits_{\tau =1}^{t-1}B_{t, \tau }^{(3)}\xi _{\tau
}-\sum\limits_{\tau =1}^{t-1}A_{t, \tau }^{(3)}z_{\tau }(\xi_{1:\tau -1})\right\}. \label{xtdef}
\end{equation}%
Here, we adopt the previous notation $z_{1:t-1}:=\left( z_{1},\ldots
,z_{t-1}\right) $ from $\xi $. By $\Pi $ we denote the operator which
maps a policy $y:=\left( y_{t}\left( \xi _{1:t-1}\right) \right)
_{t=1,\ldots ,T}$ to a new policy $z:=\Pi (y)$ defined iteratively 
by 
\begin{equation}
z_{t}\left( \xi _{1:t-1}\right) :=\pi _{X_{t}\left( z_{1:t-1},\xi
_{1:t-1}\right) }\left( y_{t}\left( \xi _{1:t-1}\right) \right) \quad
\forall \xi, \,\,\forall t=1,\ldots ,T,  \label{projdef}
\end{equation}
starting from $z_{1} :=\pi _{X_{1}}\left( y_{1}\right)$.
For example, for $t=1,2,3,\ldots $ one gets successively that%
\begin{eqnarray*}
z_{1} &:&=\pi _{X_{1}}\left( y_{1}\right) , \\
z_{2}\left( \xi _{1}\right) &:&=\pi _{X_{2}\left( z_{1},\xi _{1}\right)
}\left( y_{2}\left( \xi _{1}\right) \right), \quad \forall \xi _{1}, \\
z_{3}\left( \xi _{1},\xi _{2}\right) &:&=\pi _{X_{3}\left( z_{1},z_{2}\left(
\xi _{1}\right) ,\xi _{1},\xi _{2}\right) }\left( y_{3}\left( \xi _{1},\xi
_{2}\right) \right), \quad \forall \xi _{1}\,\,\forall \xi _{2},
\end{eqnarray*}%
so that $\Pi (y)$ is correctly defined and by (\ref{xtdef}) satisfies the
hard (almost sure) constraints of (\ref{finopt}).
(\ref{projdef}) amounts to a scenario-wise projection onto the polyhedra (\ref%
{xtdef}) which can be carried out numerically by solving a convex quadratic
program subject to linear constraints. In the special case of rectangular
sets $[y^{lo},y^{up}]$, which can be modeled as a hard constraint in \eqref{m1def} by putting for $t=1,\ldots ,T$ and $\tau=1,\ldots ,t-1$:
\begin{equation}\label{boxexp}
A_{t, t}^{(3)}:=(I,-I)^{T},\,b_t^{(3)}:= \left( \begin{array}{r}y^{up}_t\\ -y^{lo}_t \end{array} \right),\,A_{t,\tau}^{(3)}:=0,\,B_{t,\tau}^{(3)}:=0,
\end{equation}
an explicit formula can be exploited: projection of a policy then just means cutting it off at
the given lower and upper limits. For instance, in the context of the hard
constraints in (\ref{hydrototal}), one has that%
\begin{equation}
\Pi (y)=\Pi (y_{1},y_{2}\left( \cdot \right) )=\left( \max \{y^{lo},\min
\{y_{1},y^{up}\}\},\max \{y^{lo},\min \{y_{2}\left( \cdot \right)
,y^{up}\}\}\right) .  \label{projhydro}
\end{equation}%
As mentioned above, projection via $\Pi $ is a way to enforce the hard
constraints. This offers several alternatives to the above-mentioned direct intersection
of feasible policies from $M_1$ with a given (typically finite-dimensional) subclass
$\mathcal{K}$. One option would consist in working from the very beginning with projected policies so that the feasible set would become  $M_1\cap\Pi (\mathcal{K})$ rather than $M_1\cap\mathcal{K}$. Indeed, we shall see in Lemma \ref{subpol} that the intersection with the original 
infinite-dimensional feasible set may be substantially larger by doing so (in particular it would be no more empty in the example discussed before). A second option would consist in relaxing the hard constraints to probabilistic constraints similar to the ones given from the beginning and projecting them afterwards onto the set defined by hard constraints. We formalize this idea by introducing the alternative 
(infinite-dimensional) constraint set 
\begin{align}
& M_{2}:=\Big\{\left( y_{t}\left( \xi _{1:t-1}\right) \right) _{t=1,\ldots ,T}|
& &  \label{m2def} \\
& \mathbb{P}\left( \left. 
\begin{array}{r}
\sum\limits_{\tau =1}^{t}A_{t, \tau }^{(2)}y_{\tau }\left( \xi _{1:\tau
-1}\right) +\sum\limits_{\tau =1}^{t}B_{t, \tau }^{(2)}\xi _{\tau }\leq b_{t}^{(2)}
\\ 
\sum\limits_{\tau =1}^{t}A_{t, \tau }^{(3)}y_{\tau }\left( \xi _{1:\tau
-1}\right) +\sum\limits_{\tau =1}^{t-1}B_{t, \tau }^{(3)}\xi _{\tau }\leq b_{t}^{(3)}%
\end{array}%
\right\} t=1,\ldots ,T\right) \geq p\quad \Big\}. & &  \notag
\end{align}%
We shall see in Lemma \ref{projlem} that the projection of $M_{2}$ onto the hard constraints yields the set $M_{1}$, so there is no difference in the solution of \eqref{finopt} in the original infinite-dimensional setting. When considering intersections with a subclass $\mathcal{K}$, however, a significant advantage over working with 
$M_1$ may be observed. 
\subsection{Approximating the original problem by means of subclasses of decision rules}
The following result clarifies the relations between the feasible sets $M_1$, $M_2$ introduced above and their intersection with (projections of) subclasses of decision rules:
\begin{lem}\label{subpol}
\label{projlem} If $\mathcal{K}$ is an
arbitrary subset of Borel measurable policies $\left( y_{t}\left( \xi
_{1:t-1}\right) \right) _{t=1,\ldots ,T}$, then the following chain of
inclusions holds true: 
\begin{equation*}
M_{1}\cap \mathcal{K}\subseteq \Pi (M_{2}\cap \mathcal{K})\subseteq
M_{1}\cap \Pi (\mathcal{K})\subseteq M_{1}.
\end{equation*}%
In particular, by setting $\mathcal{K}$ equal to the space of all Borel
measurable policies, we derive that $\Pi (M_{2})=M_{1}$.
\end{lem}

{\textbf{Proof.}} Let $z\in M_{1}\cap \mathcal{K}$. Then, the probabilistic constraint for the
first and the almost sure constraints for the other inequality system in (\ref%
{m1def}), respectively, guarantee that the joint probabilistic constraint in
(\ref{m2def}) is satisfied, hence $z\in M_{2}\cap \mathcal{K}$. With $z$
fulfilling the almost sure constraints in (\ref{m1def}), we have that $z=\Pi
(z)$, whence $z\in \Pi (M_{2}\cap \mathcal{K})$. This proves the first
inclusion in the above chain. Next, as for the second inequality, let $z\in
\Pi (M_{2}\cap \mathcal{K})$, hence $z=\Pi (y)$ for some $y\in M_{2}\cap 
\mathcal{K}$. In particular, $z\in \Pi (\mathcal{K})$ and it remains to show
that $z\in M_{1}$. As an image of the mapping $\Pi $, $z$ satisfies the
almost sure constraints of (\ref{m1def}). By $y\in M_{2}$ and (\ref{m2def}),
there exists a measurable set $S\subseteq \Omega $ such that $\mathbb{P}%
\left( S \right) \geq p$ and 
\begin{eqnarray*}
\sum\limits_{\tau =1}^{t}A_{t, \tau }^{(2)}y_{\tau }\left( \xi _{1:\tau
-1}\left( \omega \right) \right) +\sum\limits_{\tau =1}^{t}B_{t, \tau
}^{(2)}\xi _{\tau }\left( \omega \right) &\leq& b_{t}^{(2)}\\
\sum\limits_{\tau
=1}^{t}A_{t, \tau }^{(3)}y_{\tau }\left( \xi _{1:\tau -1}\left( \omega \right)
\right) +\sum\limits_{\tau =1}^{t-1}B_{t, \tau }^{(3)}\xi _{\tau }\left(
\omega \right) &\leq& b_{t}^{(3)}
\end{eqnarray*}%
are satisfied for all $t=1,\ldots ,T$ and all $\omega \in S$. By (\ref%
{projdef}), the second inequality system implies (successively for $t$ from $%
1$ to $T$) that%
\begin{equation*}
y_{t}\left( \xi _{1:t-1}\left( \omega \right) \right) \in X_{t}\left(
z_{1:t-1},\xi _{1:t-1}\left( \omega \right) \right) \quad \forall t=1,\ldots
,T,\,\,\forall \omega \in S.
\end{equation*}%
Hence, again by (\ref{projdef}), $\left( z_{t}\left( \xi _{1:t-1}\right)
\left( \omega \right) \right) _{t=1,\ldots ,T}=\left( y_{t}\left( \xi
_{1:t-1}\right) \left( \omega \right) \right) _{t=1,\ldots ,T}$ for all $%
\omega \in S$. Therefore, the first inequality system above can be
written as%
\begin{equation*}
\sum\limits_{\tau =1}^{t}A_{t, \tau }^{(2)}z_{\tau }\left( \xi _{1:\tau
-1}\left( \omega \right) \right) +\sum\limits_{\tau =1}^{t}B_{t, \tau
}^{(2)}\xi _{\tau }\left( \omega \right) \leq b_{t}^{(2)}\quad \forall t=1,\ldots
,T, \,\,\forall \omega \in S.
\end{equation*}%
Since $\mathbb{P}\left(S \right) \geq p$ it follows that $z$ satisfies
the probabilistic constraint in (\ref{m1def}). Summarizing we have shown
that also $z\in M_{1}$, whence the desired inclusion follows. The last
inclusion is trivial.\hfill $\square$

\bigskip\noindent 
The previous lemma suggests to consider the following four
optimization problems each of them being some relaxation of our original
optimization problem (\ref{finopt}):%
\begin{eqnarray}
&&\min \{h(y)|y \in M_{1}\cap \mathcal{K}\},  \label{opt1} \\
&&\min \{h(z)|z \in \Pi (\arg \min \{h(y)|y\in M_{2}\cap \mathcal{K}\})\},
\label{opt2} \\
&&\min \{h(y)|y\in \Pi (M_{2}\cap \mathcal{K})\},  \label{opt3} \\
&&\min \{h(y)|y \in M_{1}\cap \Pi (\mathcal{K})\}.  \label{opt4}
\end{eqnarray}

\noindent Here $h$ refers to the objective function of (\ref{finopt}) and $%
\mathcal{K}$ is a given subclass of decision policies. The meaning of (\ref%
{opt1}), (\ref{opt3}) and (\ref{opt4}) is clear and relates to the feasible
sets considered in Lemma \ref{projlem}. In (\ref{opt2}) we determine first
the solution(s) of the inner optimization problem $\min \{h(y)|y\in
M_{2}\cap \mathcal{K}\}$ and then project them via $\Pi $. If this inner
optimization problem has multiple solutions, then we choose those of their
projections under $\Pi $ yielding the smallest value of the objective. We
observe that (\ref{opt3}) has the same optimal value as the problem 
\begin{equation}
\min \{h(\Pi (y))|y\in M_{2}\cap \mathcal{K}\},  \label{opt5}
\end{equation}%
where the projection is shifted from the constraints to the objective, and
that $y$ is a solution of (\ref{opt5}) if and only if $\Pi (y)$ is \ a
solution of (\ref{opt3}). Hence, (\ref{opt3}) and (\ref{opt5}) are
equivalent and it may be a matter of convenience which of the two forms is
preferred. The potential advantage of (\ref{opt2}) say over (\ref{opt3}) and
(\ref{opt4}) is that projections do not have to be dealt with in the
constraints or in the objective directly but can be carried out after
solving the problem.

\begin{lem}
\label{optvalchain}Denote by $\varphi _{1}$, $\varphi _{2}$, $\varphi _{3}$, 
$\varphi _{4}$, respectively, the optimal values of problems (\ref{opt1})-(%
\ref{opt4}) and by $\varphi $ the optimal value of the originally given
problem (\ref{finopt}). Then, any solution of problems (\ref{opt1})-(\ref%
{opt4}) is feasible for problem (\ref{finopt}) and it holds that 
\begin{equation*}
\varphi _{1},\varphi _{2}\geq \varphi _{3}\geq \varphi _{4}\geq \varphi .
\end{equation*}
\end{lem}

{\textbf{Proof.}} From Lemma \ref{projlem} we see that any feasible point and, hence, any
solution of (\ref{opt1}), (\ref{opt3}) and (\ref{opt4}) is feasible for (\ref%
{finopt}). From the inclusions of Lemma \ref{projlem} it follows that $%
\varphi _{1}\geq \varphi _{3}\geq \varphi _{4}\geq \varphi $. Now, let $%
z^{\ast }$ be a solution of (\ref{opt2}). Then, there exists some $y^{\ast
}\in M_{2}\cap \mathcal{K}$ such that $z^{\ast }=\Pi \left( y^{\ast }\right) 
$ and $y^{\ast }$ solves the problem $\min \{h(y)|y\in M_{2}\cap \mathcal{K}%
\}$. In particular, $z^{\ast }\in \Pi \left( M_{2}\cap \mathcal{K}\right) $
is feasible for (\ref{opt3}). This implies first $z^{\ast }\in M_{1}$ by
Lemma \ref{projlem} and, hence, the asserted feasibility of $z^{\ast }$ for (%
\ref{finopt}). Second, it implies the desired remaining relation $\varphi
_{2}=h(z^{\ast })\geq \varphi _{3}$.\hfill $\square$

\bigskip\noindent 
Lemma \ref{optvalchain} can be interpreted as follows: Problem (%
\ref{opt1}) reflects the pure transition to a subclass $\mathcal{K}$ of
policies in the originally given problem (\ref{finopt}). The resulting loss
in optimal value equals $\varphi _{1}-\varphi \geq 0$. In contrast, using
projections onto hard constraints in the one or other way as in (\ref{opt3})
and (\ref{opt4}) may lead to smaller losses in the optimal values. Of
course, this advantage of working with projections requires that the
computational gain by passing to an interesting subclass $\mathcal{K}$ is
not destroyed by the projection procedure. This is why in Section \ref{lindec} we
shall introduce the class of linear decision rules as a suitable one
harmonizing well to a certain degree with projections onto polyhedral sets.
The following example illustrates Lemma \ref{optvalchain}:

\begin{ex}\label{ex1}
Consider the following problem with policies $y_{1},y_{2}(\xi _{1})$ as
variables: 
\begin{eqnarray*}
&&\min \; y_{1}\quad\mbox{ \rm subject to} \\
&&\mathbb{P}(\xi _{1}\leq y_{1},\,\,\xi _{2}\leq y_{2}(\xi _{1}))\geq p \\
&&y_{1},\,y_{2}(\xi _{1})\in [0, 1],\quad \mathbb{P}-%
\mbox{\rm almost
surely.}
\end{eqnarray*}%
We assume that the random vector $\xi =\left( \xi _{1},\xi _{2}\right) $
follows a uniform distribution over the set $\Theta = \left( \left[ -1,1\right]
\times \left[ 0,1\right] \right) \cup \left( \left[ 0,1\right] \times \left[
0,-1\right] \right) $ and that $p=1/3$. As a subclass of
policies, we consider (purely) linear second stage decisions:%
\begin{equation*}
\mathcal{K}:=\{\left( y_{1},y_{2}(\xi _{1})\right) |\exists a \geq -1%
:y_{2}(\xi _{1})=a\xi _{1}\}.
\end{equation*}

\begin{itemize}
\item \textbf{Solution of the original problem }(\ref{finopt}):

We claim that the optimal value $\varphi $ of the original problem equals $0$. 
Indeed, it cannot be smaller than $0$ due to the constraint $y_{1}\geq 0$.
On the other hand, $y_{1}:=0$ and $y_{2}(\xi _{1}):=1$ for all $\xi _{1}$
represents a feasible policy because it clearly satisfies the almost sure
constraints and the set of $\xi $ satisfying $\xi _{1}\leq 0$ and \,\,$\xi
_{2}\leq 1$ covers one-third of the support of $\xi $.
Hence the probabilistic constraint is satisfied too. The objective value
associated with this feasible policy equals $y_{1}=0$, so $\varphi =0$ as
asserted.

\item \textbf{Solution of problem }(\ref{opt1}):

The feasible set here is $M_{1}\cap \mathcal{K}$ and a feasible second stage
policy $y_{2}(\xi _{1})=a\xi _{1}$ has to be trivial $\left( a=0\right) $ in
order to satisfy the almost sure constraint $0\leq \,y_{2}(\xi _{1})\leq 1$.
Then, the only choice for $y_{1}$ such that $\left( y_{1},0\right) $
satisfies the probabilistic constraint is $y_{1}:=1$ (only then, the set of $%
\xi $ satisfying $\xi _{1}\leq y_{1}$ and$\,\,\xi _{2}\leq 0$ covers one-third of the support of $\xi $). Hence the feasible set in this problem
reduces to a singleton and its optimal value equals to the objective value
of this singleton: $\varphi _{1}=y_{1}=1$.

\item \textbf{Solution of problems} (\ref{opt2}) and (\ref{opt3}):

As stated above, (\ref{opt3}) is equivalent with (\ref{opt5}). In our
example, $h$ is the projection onto the first component, hence we seek to
minimize $\left( \Pi (y)\right) _{1}$ over the constraint set%
\begin{eqnarray}
M_{2}\cap \mathcal{K}&=&\{\left( y_{1},a\xi _{1}\right) \mid a \geq -1, \mathbb{P}(y_{1},a\xi _{1}\in [0,1],\,\,\xi_{1}\leq y_{1},\,\,\xi _{2}\leq a\xi _{1})\geq 1/3\}\nonumber \\
& =&  \{(y_1, a \xi_1) \mid y_1 \in [0,1], a \geq -1, \psi (y_1, a) \geq 1/3\}\label{proconopt3}
\end{eqnarray}
where $\psi (y_1, a):= \mathbb{P}((\xi_1, \xi_2) \in \mathcal{S}(a, y_1))$
with  
\[
\mathcal{S}(a, y_1):=
\left\{(\xi_1, \xi_2)\in \Theta\mid \xi_1 \leq y_1, \xi_2 \leq a \xi_1, 0 \leq a \xi_1 \leq 1 \right\}
\]
(see Figure \ref{figs1}).
%\begin{figure}[H]
%\begin{tabular}{ll}
%\input{Case12_1.pstex_t}& \input{Case14_1.pstex_t} \\
%\input{Case12_2.pstex_t}& \input{Case12_3.pstex_t} 
%\end{tabular}
%\caption{Representations of $\mathcal{S}(a, y_1)$ and ${\tilde{\mathcal{S}}}(a, y_1)$:
%top figures for $-1 \leq a \leq 0$, bottom left for $0 < a \leq 1$, and bottom right for $a>1$.}
%\label{figs1}
%\end{figure}
\begin{figure}[bt]
\begin{center}
\includegraphics[width=12cm]{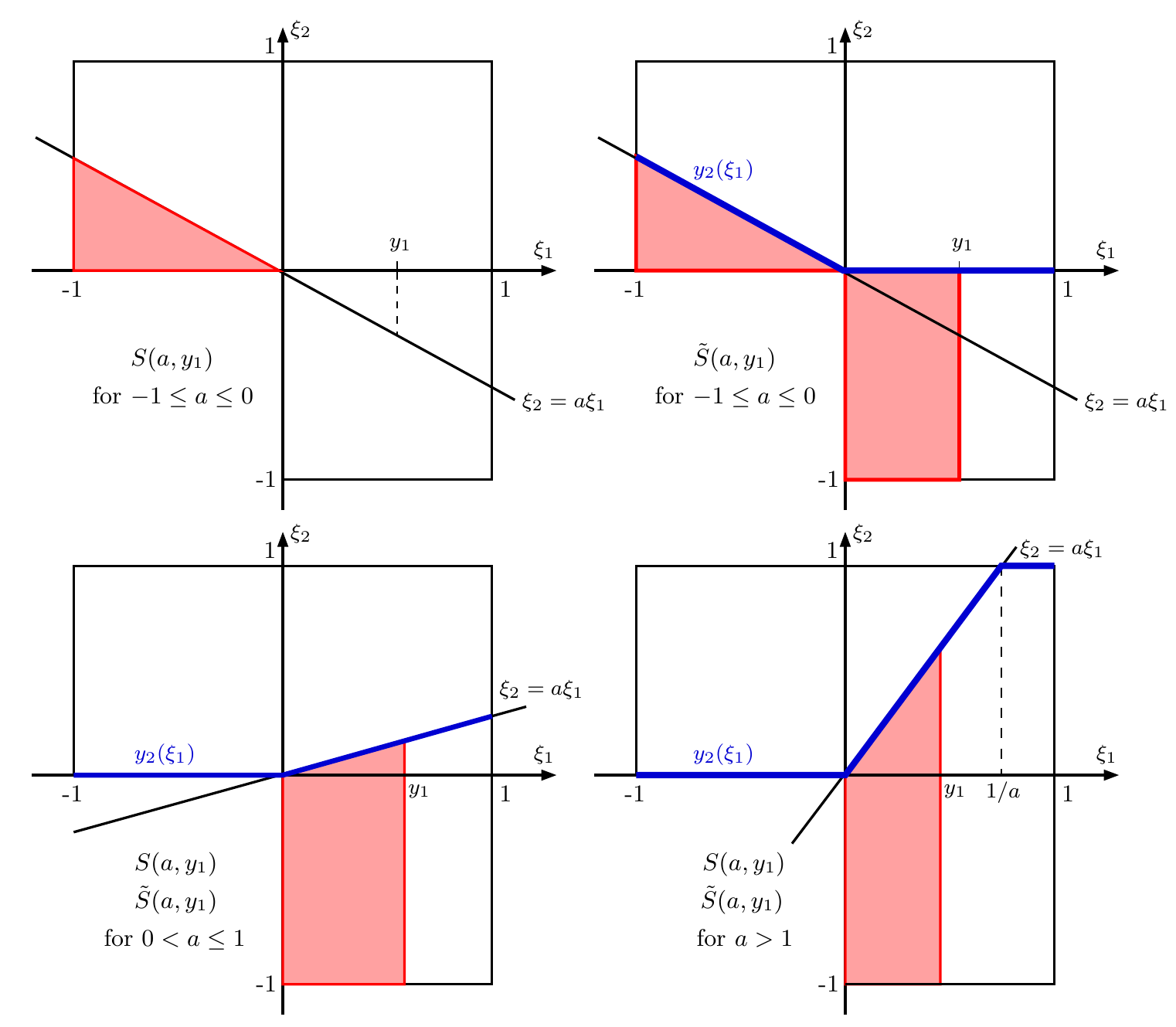}
\end{center}
\caption{Representations of $\mathcal{S}(a, y_1)$ and ${\tilde{\mathcal{S}}}(a, y_1)$:
top figures for $-1 \leq a \leq 0$, bottom left for $0 < a \leq 1$, and bottom right for $a>1$.}
\label{figs1}
\end{figure}
Note that in \eqref{proconopt3} we were allowed to extract the deterministic constraint $y_1 \in [0,1]$ from the probabilistic constraint. 

As $\left( \Pi (y)\right) _{1}$ is the projection of $%
y_{1}$ onto the first stage almost sure constraint set $X_{1}=\left[ 0,1%
\right] $ (see (\ref{projdef}) and (\ref{xtdef})), we get that $\left( \Pi
(y)\right) _{1}=y_{1}$. Consequently, according to (\ref{opt5}), we want to
minimize $y_{1}$ for all policies $\left( y_{1},a\xi _{1}\right) $ belonging
to (\ref{proconopt3}).
We consider three cases:
\begin{itemize}
\item[(i)] For $-1 \leq a \leq 0$ (see top left in Figure \ref{figs1}), we have
$\psi (y_1, a)=-a/6<1/3$. 
\item[(ii)] For $0 < a \leq 1$ (see  bottom left in Figure \ref{figs1}), we have  
$\psi (y_1, a)=\frac{1}{3}(y_1 + a y_{1}^{2}/2)$. The smallest value of $y_1$ satisfying 
$\psi (y_1, a) \geq 1/3$ is obtained taking $a=1$ and $y_1=-1+\sqrt{3}>2/3$.
\item[(iii)] For $a>1$ (see  bottom right in Figure \ref{figs1}), we get
\[
\psi (y_1, a)=\left\{
\begin{array}{ll}
\frac{1}{3}(y_1 + a y_{1}^{2}/2) & \mbox{if }y_1 \leq 1/a,\\
 \frac{1}{2 a} & \mbox{otherwise.}
\end{array}
\right.
\]
In particular, $\psi (\frac{2}{3},\frac{3}{2})=\frac{1}{3}$. We distinguish the two subcases:
\begin{itemize}
\item[(1)] $a>3/2$: if $y_1>1/a$ then $\psi (y_1, a)=\frac{1}{2 a}<\frac{1}{3}$
and if $0 \leq y_1 \leq \frac{1}{a}$ then  $\psi (y_1, a)\leq \psi (\frac{1}{a}, a) =\frac{1}{2 a}<\frac{1}{3}$.
\item[(2)] $1 \leq a \leq 3/2$: If $0\leq y_1 <\frac{2}{3}$ then $y_1\leq 1/a$ and, hence,
\[
\psi (y_1, a)<\frac{1}{3}\left(\frac{2}{3} + \frac{3}{2} \frac{4}{2\cdot 9}\right)=
\frac{1}{3}
\]
\end{itemize}
\end{itemize}
Summarizing, the best value of the objective at an admissible solution of (\ref{opt5}) equals $\frac{2}{3}$ and is realized uniquely by the optimal policy $\left( \frac{2}{3},\frac{3}{2}\xi
_{1}\right) $. The latter is therefore the unique optimal solution of (\ref{opt5}). According to our observation above, its projection 
\begin{equation}\label{projsol}
\Pi \left( \frac{2}{3},\frac{3}{2}\xi
_{1}\right)=\left( \frac{2}{3},\max \{0,\min \{\frac{3}{2}\xi _{1},1\}\}\right)
\end{equation}
onto the almost sure constraints in our example is an optimal
solution of (\ref{opt3}). The associated function value equals $h\left( \frac{2}{3},\frac{3}{2}\xi_{1}\right)=2/3$ which therefore is the optimal value of \eqref{opt3}. It follows that $\varphi_{3}=2/3$. 

On the other hand, as we have already observed that \,$h(y)=h(\Pi (y))=y_1$\,
due to \,$0\leq y_{1}\leq 1$, it follows that the unique optimal solution 
$\left( \frac{2}{3},\frac{3}{2}\xi_{1}\right) $ of (\ref{opt5}) yields the unique optimal solution to the problem
\[
\min \{h(y)|y\in M_{2}\cap \mathcal{K}\}
\]
at the same time. Hence, its projection onto the almost sure constraints is the already identified solution (\ref{projsol}) of problem (\ref{opt3}) implying that the optimal value of
problem (\ref{opt2}) is the same as that of (\ref{opt3}): $\varphi_2=\varphi
_{3}=2/3$.
\item \textbf{Solution of problem }(\ref{opt4}): By virtue of \eqref{projhydro}, the 
policies belonging to the set $\Pi (\mathcal{K})$ have the form 
$(y_{1},\max \{0,\min \{a\xi_{1},1\}\})$
for some $y_1\in [0,1]$ and $a\geq -1$ (see Figure \ref{figs1}). Since 
these policies already satisfy the almost sure
constraints, all one has to add in order to get a policy feasible for (\ref%
{opt4}) is the satisfaction of the probabilistic constraint. Observe that
\[
M_1 \cap \Pi(\mathcal{K})  =  \{(y_1,\max \{0,\min \{a\xi_{1},1\}\})\mid y_1\in [0,1],\, a \geq -1,\,\tilde{\psi}(y_1, a) \geq 1/3\} 
\]
where $\tilde{\psi}(y_1, a):= \mathbb{P}((\xi_1, \xi_2) \in {\tilde{\mathcal{S}}}(a, y_1))$
with  
\[
{\tilde{\mathcal{S}}}(a, y_1)=
\left\{(\xi_1, \xi_2) \in \Theta\mid \xi_1 \leq y_1, \xi_2 \leq \max(0, \min(a\xi_1, 1)) \right\}
\]
(see Figure \ref{figs1}). For $-1 \leq a \leq 0$ (see Figure \ref{figs1}), we have
\[
\tilde{\psi}(y_1, a)=\frac{1}{3}(y_1 - \frac{a}{2})< \tilde{\psi}(1/2, -1)=\frac{1}{3}
\quad\forall y_1<1/2. 
\]
For $0 < a \leq 1$ (see bottom left in Figure \ref{figs1}), we have  
$\tilde{\psi}(y_1, a)=\frac{1}{3}(y_1 + a y_{1}^{2}/2)$. The smallest value of $y_1$ satisfying 
$\tilde{\psi}(y_1, a) \geq 1/3$ is obtained taking $a=1$ and $y_1=-1+\sqrt{3}>1/2$.
Finally, for $a>1$ (see bottom right in Figure \ref{figs1}), we assume that $y_1 \leq 1/2$.
Then,
\begin{eqnarray*}
y_1>1/a&\Longrightarrow&\tilde{\psi}(y_1, a)=\frac{1}{3}(2y_1-\frac{1}{2a})< \frac{y_1}{3} \leq \frac{1}{3}\\
y_1 \leq 1/a&\Longrightarrow&\tilde{\psi}(y_1, a)=\frac{1}{3}(y_1 +\frac{a y_1^2}{2}) \leq \frac{1}{3}(\frac{1}{2} + \frac{1}{2a})<\frac{1}{3}.
\end{eqnarray*}
This means that there is no feasible policy with $y_1 \leq 1/2$ and $a>1$.
Consequently, the optimal value of (\ref{opt4}) equals $\varphi _{4}=1/2$ and is realized by the policy $(1/2,\max \{0,\min \{-\xi_{1},1\}\})$ which is the projection of the  decision rule $(1/2,-\xi_{1})\in\mathcal{K}$ onto the hard box constraints.
\end{itemize}
\end{ex}

\section{Probabilistic model and linear decision rules}\label{probmod}

Example \ref{ex1} has illustrated the different approximating optimization problems with respect to the given one \eqref{finopt}. In order to formulate these ideas in a practically
meaningful framework, one has to specify the probabilistic model for the random vector $\xi$ and a suitable subclass $\mathcal{K}$ of decision rules in Lemma \ref{subpol}.
\subsection{Probabilistic model} \label{gensetting}
We introduce in this section the class of stochastic processes $(\xi_t)$ we consider. Each component $\xi_t(m)$ of $\xi_t$ follows a linear model of the form
\begin{equation} \label{modeletimeseries}
\begin{array}{l}
\displaystyle{\sum_{k=0}^{p_t(m)}} \; \alpha_{t, k}(m) \xi_{t-k}(m) =\mu_t(m)+ \displaystyle{\sum_{k=0}^{q_t(m)}} \; \beta_{t, k}(m) \varepsilon_{t-k}(m),\;m=1,\ldots,M,
\end{array}
\end{equation}
where $\mu_t$ is the tendency for period $t$ and lags $p_t(m), q_t(m)$ are nonnegative and depend on time.
We assume that for every $t$, the coefficients  $\alpha_{t, 0}(m), \alpha_{t, p_t(m)}(m)$, and $\beta_{t, q_t(m)}(m)$ are nonzero. 

Finally, the noises are supposed to obey centered Gaussian laws $\varepsilon_t\sim\mathcal{N}(0,\Sigma_t)$, pairwise independent for different time steps. We recall the notation $\mathcal{N}(\mu ,\Sigma)$ for referring to a multivariate Gaussian distribution with mean $\mu$ and covariance matrix $\Sigma$. Hence, $\varepsilon:=(\varepsilon_1,\ldots ,\varepsilon_T)\sim\mathcal{N}(0,\Sigma)$, where $\Sigma$ is a block-diagonal covariance matrix whose blocks are the covariance matrices $\Sigma_t$ of the components $\varepsilon_t$.
\begin{rem} We assume that the parameters of model \eqref{modeletimeseries} are known.
In its full generality, model \eqref{modeletimeseries} is not identifiable.
Additional assumptions are needed to identify lags $p_t(m), q_t(m)$ and calibrate the model parameters.
As special cases, the identifiable SARIMA (with constant lags) and 
Periodic Autoregressive (PAR, with periodic time-dependent lags) models can 
be considered.
\end{rem}
Using iteratively model equation \eqref{modeletimeseries}, for each instant $t=1,\ldots,T$, we can decompose $\xi_t(m)$ as a function
of noises $\varepsilon_1,\ldots,\varepsilon_t$ and of
past observations of the process $(\xi_t)$ and of the noises (observations for instants $0,-1,-2,\ldots$). More precisely, for 
every $t=1,\ldots ,T$ and for every component $m$, we have for $\xi_t(m)$ a decomposition of the form
\begin{equation} \label{eqxtm}
\xi_t(m)=c_t(m)+\displaystyle{\sum_{k=1}^{r_t(m)}} \;\gamma_{t, k}(m) \xi_{1-k}(m)+\displaystyle{\sum_{k=1}^{s_t(m)}} \;\delta_{t, k}(m) \varepsilon_{1-k}(m)+\displaystyle{\sum_{k=1}^{t}} \;\theta_{t, k}(m) \varepsilon_{k}(m)
\end{equation}
for some lags $r_t(m)$ and $s_t(m)$ that represent the minimal number of past observations of 
respectively the stochastic processes $(\xi_t)$ and $(\varepsilon_t)$ that are necessary
to decompose $\xi_t(m)$ over its past.
This decomposition will be used in the next sections. In this decomposition, the first two sums gather the past realizations of process $(\xi_t)$ and of the noises.
Lemma \ref{lemcalcijh} stated and proved in the Appendix, provides the formulae to compute iteratively the coefficients appearing in the decompositions
of $\xi_1(m), \xi_2(m)$, $\ldots,\xi_T(m), m=1,\ldots,M$, of the form \eqref{eqxtm} above. The computation of these
coefficients is necessary when one is interested in solving the optimization problems we consider in the next sections
when $(\xi_t)$ is of the form \eqref{modeletimeseries}. A similar decomposition for less general models was given 
in \cite{gsorl2012}, \cite{gs2012}.

It is convenient to write \eqref{eqxtm} in the compact form 
\begin{equation} \label{substi2}
\xi_t={\tilde \mu_t}+ \Theta_t \varepsilon\quad (t=1,\ldots ,T),
\end{equation}
where for each $t=1,\ldots ,T$,
\begin{itemize}
\item ${\tilde \mu_t}$ is a constant vector in $\mathbb{R}^M$ with component $m$ given by
$$
{\tilde \mu_t}(m)=c_t(m)+\displaystyle{\sum_{k=1}^{r_t(m)}} \;\gamma_{t, k}(m) \xi_{1-k}(m)+\displaystyle{\sum_{k=1}^{s_t(m)}} \;\delta_{t, k}(m) \varepsilon_{1-k}(m),
$$
\item $\Theta_t$ is the $M\small{\times}M T$ matrix
$$
\Theta_t=\left(\mbox{diag}(\theta_{t, 1}(1), \ldots, \theta_{t, 1}(M)), \ldots, \mbox{diag}(\theta_{t, t}(1),\ldots,\theta_{t, t}(M)), 0_{M \small{\times}M(T-t)}\right)
$$
where the coefficients $\theta_{t, j}(m)$ are given in Lemma \ref{lemcalcijh}.
\end{itemize}
\subsection{Linear decision rules}\label{lindec}
As mentioned in Section \ref{hardwait} the numerical solution of problem \eqref{finopt}
requires to reduce the space of all Borel measurable decision policies to some convenient finite-dimensional subspace. A simple and widely used way to do so consists in considering so-called {\it linear decision rules} as policies which are defined as the set
\begin{equation}\label{ldr}
\mathcal{K}:=\{\left( y_{t}\left( \xi _{1:t-1}\right) \right) _{t=1,\ldots ,T}\mid
\exists F_t,f_t:\,y_{t}\left( \xi _{1:t-1}\right)=F_t\xi _{1:t-1}+f_t\quad (t=1,\ldots ,T)\},
\end{equation}
with matrices $F_t$ and vectors $f_t$ of appropriate size. Since the first stage decision $y_1$ is deterministic, we convene about fixing $F_1:=0$.
\subsubsection{The random inequality system under linear decision rules}
Under linear decision rules and the probabilistic model \eqref{substi2}, our generic random inequality system
\begin{equation}\label{generic}
\sum\limits_{\tau =1}^{t}A_{t, \tau }y_{\tau }\left( \xi _{1:\tau -1}\right)
+\sum\limits_{\tau =1}^{t}B_{t, \tau }\xi _{\tau }\leq b_{t}\quad t=1,\ldots
,T 
\end{equation}
turns into (for $t=1,\ldots ,T$) 
\begin{equation}\label{randlin}
 \underbrace{\left(\sum_{\tau=1}^t \; A_{t, \tau} F_{\tau} \Theta_{1:\tau-1}+B_{t, \tau} \Theta_{\tau} \right)}_{G_t(x)}\varepsilon \leq \underbrace{b_t-\sum_{\tau=1}^t B_{t, \tau} {\tilde \mu}_{\tau}-\sum_{\tau=1}^t A_{t, \tau} f_{\tau} - \sum_{\tau=1}^t \; A_{t, \tau} F_{\tau} {\tilde \mu}_{1:\tau-1}}_{g_t(x)}.
\end{equation}
In this system, $\varepsilon$ is the transformed random vector, whereas now $x:=(F_t,f_t)_{t=1,\ldots T}$ 
represents a finite-dimensional decision vector approximating the original decision policies $\left( y_{t}\left( \xi _{1:t-1}\right) \right) _{t=1,\ldots ,T}$. With the notation introduced below the corresponding expressions, we may compactly rewrite \eqref{randlin} in the form
\begin{equation}\label{compact}
G_t(x)\varepsilon\leq g_t(x)\quad (t=1,\ldots ,T),
\end{equation}
where the $G_t, g_t$ are affine linear mappings of $x$. When relating these mappings not to the generic system \eqref{generic} but to the concrete systems of hard and soft constraints in \eqref{finopt} labeled by upper indices (1), (2), (3), we shall use the corresponding upper indices for the mappings  
$G_t$ and $h_t$ as well.

We observe that thanks to affine linearity of $G_t, g_t$, the set of $x$ satisfying \eqref{compact} is convex for each fixed $\varepsilon$.
\subsubsection{The objective function under linear decision rules}
From \eqref{substi2} and $\varepsilon$ having a centered distribution, it follows that the expectation of $\xi_t$ equals $\tilde{\mu}_t$. Therefore,
the objective of our problem \eqref{finopt} takes under linear decision rules the form
\[
\underbrace{\sum_{t=1}^T\langle h_t,F_t\tilde{\mu}_{1:t-1}+f_t\rangle}_{\mathcal{J}_1(x)} + 
\sum_{t=1}^T\left\langle\mathcal{P}_t,\mathbb{E}\left( \sum\limits_{\tau
=1}^{t}A_{t, \tau }^{(1)} y_{\tau }(\xi _{1:\tau -1})+\sum\limits_{\tau
=1}^{t}B_{t, \tau }^{(1)}\xi _{\tau }-b_{t}^{(1)}\right) _{+}\right\rangle
\]
where in the definition of $\mathcal{J}_1$ we used once more the convention $x:=(F_t,f_t)_{t=1,\ldots T}$. Now, applying \eqref{compact} with upper index (1) referring to the inequality subsystem penalized in the objective, we can rewrite the objective of \eqref{finopt} under linear decision rules as 
$\mathcal{J} (x):=\mathcal{J}_1(x)+\mathcal{J}_2(x)$, where
\[
\mathcal{J}_2(x):=\sum_{t=1}^T\left\langle\mathcal{P}_t,\mathbb{E}\left( G_t^{(1)}(x)\varepsilon-  g_t^{(1)}(x)\right) _{+}\right\rangle
\] 
\begin{lem}\label{betconv}
$\mathcal{J}$ is convex. 
\end{lem}
{\textbf{Proof.}}
Since $\mathcal{J}_1$ is linear, it suffices to check convexity of $\mathcal{J}_2$. As mentioned earlier,
the mappings $G_t^{(1)}, g_t^{(1)}$ are affine linear, whence the mapping 
$G_t^{(1)}(x)\varepsilon-g_t^{(1)}(x)$ is affine linear in $x$. In particular, each component of this mapping is convex in $x$ which remains true upon passing to its
maximum with zero. It follows that the components of $\mathbb{E}\left( G_t^{(1)}(x)\varepsilon-g_t^{(1)}(x)\right) _{+}$ (depending only on $x$) are convex. Now, the result follows from $\mathcal{P}_t\geq 0$.\hfill
$\square$

\bigskip\noindent
For implementation purposes, it is useful to have an analytic expression of the objective function. For this purpose, we need the folloming lemma:
\begin{lem}\label{lemmeanpiece} Let $X$ be a one-dimensional Gaussian random variable distributed according to $\mathcal{N}(m,\sigma^2)$ and let 
$a,b \in {\bar {\mathbb{R}}}$ with
$a \leq b$. Then, with $\Phi$ referring to the one-dimensional standard normal distribution function, it holds that
\begin{eqnarray*}
\mathbb{E}[\max\{a,\min\{X,b\}\}]&=&
\frac{\sigma}{\sqrt{2\pi}}\left(\exp\left(-\frac{(a-m)^2}{2\sigma^2}\right)-
\exp\left(-\frac{(b-m)^2}{2\sigma^2}\right)\right)+\\&&
(a-m)\Phi (\frac{a-m}{\sigma})+(m-b)\Phi (\frac{b-m}{\sigma})+b.
\end{eqnarray*}
\end{lem}
{\textbf{Proof.}}
With $f_X(x)=\frac{1}{\sqrt{2 \pi} \sigma} \exp \left({-\frac{(x-m)^2}{2 \sigma^2}}\right)$ being the density of $X$ and with $\tilde{\Phi}$ being the associated cumulative distribution function, we have
\begin{eqnarray*}
\mathbb{E}[\max\{a,\min\{X,b\}\}]&=&\\ 
\int_{-\infty}^aaf_X(x)dx+\int_{a}^bxf_X(x)dx+\int_b^\infty bf_X(x)dx&=&\\
a\tilde{\Phi} (a)+\int_{a}^b(x-m)f_X(x)dx+m\int_{a}^bf_X(x)dx+b(1-\tilde{\Phi} (b))&=&\\
a\tilde{\Phi} (a)+\left[\frac{-\sigma}{\sqrt{2\pi}}\exp\left(-\frac{(x-m)^2}{2\sigma^2}\right)\right]^b_a+m(\tilde{\Phi}(b)-\tilde{\Phi}(a))+b(1-\tilde{\Phi} (b))&=&\\
\frac{\sigma}{\sqrt{2\pi}}\left(\exp\left(-\frac{(a-m)^2}{2\sigma^2}\right)-\exp\left(-\frac{(b-m)^2}{2\sigma^2}\right)\right)
+(a-m)\tilde{\Phi}(a)+(m-b)\tilde{\Phi}(b)+b.
\end{eqnarray*}
On the other hand, since $\sigma^{-1}(X-m)\sim\mathcal{N}(0,1)$, we have that, for all $z$,
\[
\tilde{\Phi}(z)=\mathbb{P}(X\leq z)=\mathbb{P}(\sigma^{-1}(X-m)\leq \sigma^{-1}(z-m))=
\Phi (\sigma^{-1}(z-m))
\]
and the result follows.\hfill  $\square$

\bigskip\noindent
The only non-explicit part in our objective function $\mathcal{J}(x)$ is the vector of expectations in the definition of $\mathcal{J}_2(x)$. Its 
$i^{th}$ component is given by  
\[
\mathbb{E}[\max(X(x),0)]=\mathbb{E}[\max\{0,\min\{X(x),+\infty\}];\quad X(x):=\left( G_t^{(1)}(x)\varepsilon-g_t^{(1)}(x)\right)_i. 
\]
According to the transformation rules of Gaussian distributions, we know that 
\[
X(x)\sim\mathcal{N}(m,\sigma^2);\quad m:=-(g_t^{(1)}(x))_i;\quad
\sigma:=\sqrt{\left( G_t^{(1)}(x)\Sigma [G_t^{(1)}(x)]^T\right)_{ii}},
\]
where $\Sigma$ is the block-diagonal covariance matrix of $\varepsilon$ (see
Section \ref{probmod}). With these data, Lemma \ref{lemmeanpiece} can be employed 
(with $a:=0, b:=+\infty$) to make the objective $\mathcal{J}(x)$ fully explicit in terms of the initial data of the problem.
\subsubsection{Projection of linear decision rules onto hard constraints}\label{projldr}
The solution of problems \eqref{opt2}, \eqref{opt3}, \eqref{opt4}, \eqref{opt5} is intimately related to the ability to either explicitly or numerically compute projections $\Pi (y)$ of policies $y\in\mathcal{K}$ according to \eqref{projdef}. In the 
case of linear decision rules introduced in \eqref{ldr}, the projected policy $z:=\Pi (y)$ is obtained for 
$y=(F_t\xi _{1:t-1}+f_t)_{t=1,\ldots ,T}$ as the successive (unique) solution of  (scenario-dependent) quadratic 
optimization problems:
\begin{eqnarray}\label{quadropt}
  z_t(\xi_{1:t-1})   = & 
\left\{
\begin{array}{l}
\underset{u}{\rm argmin} \|F_t\xi _{1:t-1}+f_t-u\|^2 \\
A_{t, t}^{(3)}u\leq b_{t}^{(3)}-\sum\limits_{\tau =1}^{t-1}B_{t, \tau }^{(3)}\xi _{\tau
}-\sum\limits_{\tau =1}^{t-1}A_{t, \tau }^{(3)}z_{\tau }(\xi_{1:\tau -1}), 
\end{array}
\right.
&
\\&&\nonumber\\
&\forall \xi,\quad\forall t=1,\ldots ,T.&\nonumber
\end{eqnarray}
Here, starting from $t=1$, previously obtained solutions for $z_\tau$ are plugged in on the right-hand side of \eqref{quadropt}. Hence, for instance the first two components of $z$ are obtained as
\begin{eqnarray*}
&z_1=\underset{u}{\rm argmin}\left\{\|f_1-u\|^2|A_{1, 1}^{(3)}u\leq b_{1}^{(3)}\right\}&
\\
&z_2(\xi_1)=\underset{u}{\rm argmin}\left\{\|F_2\xi _{1}+f_2-u\|^2|A_{2, 2}^{(3)}u\leq
b_{2}^{(3)}-B_{2,1}^{(3)}\xi _{1}-A_{2,1}^{(3)}z_1\right\}\quad\forall \xi_1.&
\end{eqnarray*}
In the special case of box constraints 
\begin{equation}\label{boxhard}
y_t(\xi_{1:t-1})\in [\underline{y}_t,\overline{y}_t]\quad\mathbb{P}\mbox{-almost surely}\,\,t=1,\ldots ,T,
\end{equation}
an explicit formula for the projection of 
$y=(F_t\xi _{1:t-1}+f_t)_{t=1,\ldots ,T}$ can be provided:
\begin{equation}\label{boxform}
\Pi (y)=\left(\max\{(\underline{y}_t)_i,\min\left\{(F_t\xi _{1:t-1}+f_t)_i,(\overline{y}_t)_i\right\}\}\right)_{t=1,\ldots ,T;\,i=1,\ldots ,n_t}.
\end{equation}
\subsubsection{Probabilistic constraints under linear decision rules and Gaussian distribution}\label{proconldr}
Under the assumption of linear decision rules \eqref{ldr}, the originally dynamic probabilistic constraint
\[
\mathbb{P}\left(\sum\limits_{\tau =1}^{t}A_{t, \tau }y_{\tau }\left( \xi _{1:\tau -1}\right)
+\sum\limits_{\tau =1}^{t}B_{t, \tau }\xi _{\tau }\leq b_{t}\quad t=1,\ldots
,T \right)\geq p
\]
associated with \eqref{generic} and occurring in problems \eqref{finopt} 
turns into a conventional static probabilistic constraint
\begin{equation}\label{pcfindim}
\mathbb{P}\left(G_t(x)\varepsilon\leq g_t(x)\quad (t=1,\ldots ,T)\right)\geq p,
\end{equation}
with finite-dimensional decisions $x:=(F_t,f_t)_{t=1,\ldots T}$. \eqref{pcfindim} represents a joint linear probabilistic constraint under Gaussian distribution. This class has been intensively studied with respect to its analytical properties and numerical solution approaches, see, e.g., \cite{shap-rusc-book,prek}.
For an algorithmic treatment of such probabilistic constraints within the framework of nonlinear optimization it is important to have required information about the probability function
\[
\varphi(x):=\mathbb{P}\left(G_t(x)\varepsilon\leq g_t(x)\quad (t=1,\ldots ,T)\right)
\]
defining the inequality constraint $\varphi (x)\geq p$ in \eqref{pcfindim}. In particular, procedures computing or, better, approximating values and gradients of $\varphi$ are needed. As shown in \cite{ackooij2}, both tasks can be realized simultaneously by reduction to the computation of multivariate Gaussian distribution functions. The latter can be quite efficiently done using 
Genz' code as described in \cite{genz}. An alternative approach consists in the use of the so-called spheric-radial decomposition of Gaussian random vectors \cite{deak,royset,ackooij3}.
Another important property for algorithmic purposes is convexity of the feasible set described by \eqref{pcfindim}. While this is well known to be true in case of 
constant matrices $G_t$ and mappings $g_t$ having concave components \cite[Theorem 10.2.1]{prek}, the same does not hold true in general for \eqref{pcfindim}, in particular not for arbitrary probability levels $p$. Apart from special cases, such as the presence of one single random inequality in the system \cite{kat,panne} or specially structured covariance matrices \cite{prek2,hen}, where convexity for sufficiently large $p$ could be guaranteed, no general result on this issue seems to be available so far. 
\section{Approximating optimization problems under linear decision rules and Gaussian and truncated Gaussian distribution}
\subsection{First optimization problem}
The first optimization problem we address is \eqref{opt1}, i.e., the original problem \eqref{finopt} but with the feasible set intersected with the class of linear decision rules \eqref{ldr}. Making recourse to the compact notation introduced in Section \ref{lindec},
Problem \eqref{opt1} writes
\begin{eqnarray}
\min\{\mathcal{J}(x)&\mid&\mathbb{P}(G^{(2)}_t(x)\varepsilon\leq g^{(2)}_t(x)\quad (t=1,\ldots ,T))\geq p,\label{firstopt}\\&&G^{(3)}_t(x)\varepsilon\leq g^{(3)}_t(x)\quad (t=1,\ldots ,T),\quad\mathbb{P}\mbox{-almost surely}\}.\nonumber
\end{eqnarray}
In the definition of $g_t^{(3)}$ according to \eqref{randlin} we have to recall that
$B_{t, t}^{(3)}=0$ for all $t=1,\ldots ,T$ according to our wait-and-see perspective on hard constraints (see Section \ref{hardwait}).
\eqref{firstopt} is a nonlinear optimization problem with a joint probabilistic and a (linear) semi-infinite constraint ($\mathbb{P}$-almost surely could be replaced by 'for $\mathbb{P}$-almost all $\varepsilon\in\Xi$', where $\Xi$ is the support of the random vector $\varepsilon$). 
\begin{prop}\label{hardexp}
The hard constraint in problem \eqref{firstopt} can be explicitly represented in terms of the original data (see \eqref{randlin}) as the system of linear 
(in-)equalities for $t=1,\ldots ,T$:
\begin{eqnarray*}
\sum_{\tau=1}^t \; \left(A_{t, \tau}^{(3)} F_{\tau} \Theta_{1:\tau-1}+B_{t, \tau}^{(3)} \Theta_{\tau}\right) &=&0,\\
\sum_{\tau=1}^{t-1} B_{t, \tau}^{(3)} {\tilde \mu}_{\tau} +   \sum_{\tau=1}^t A_{t, \tau}^{(3)} f_{\tau} + \sum_{\tau=1}^t \; A_{t, \tau}^{(3)} F_{\tau} {\tilde \mu}_{1:\tau-1}&\leq&b_t^{(3)}.
\end{eqnarray*}
\end{prop}
{\textbf{Proof.}}
As mentioned above, the hard constraint in problem \eqref{firstopt} can be replaced by
\begin{equation}\label{semiinf}
G^{(3)}_t(x)\varepsilon\leq g^{(3)}_t(x)\quad\mbox{for }\mathbb{P}\mbox{-almost all }\varepsilon\in\Xi\quad (t=1,\ldots ,T).
\end{equation}
Since $\varepsilon$ follows a multivariate Gaussian distribution, its support is the whole space. 
As a consequence, some $x$ can be feasible for \eqref{semiinf} only if $G^{(3)}_t(x)=0$ which in turn implies that $g^{(3)}_t(x)\geq 0$. Conversely, any $x$ satisfying these two relations is feasible for \eqref{semiinf}. Thus, we have shown that \eqref{semiinf} is equivalent 
with the system $G^{(3)}_t(x)=0,\,g^{(3)}_t(x)\geq 0$. Now, \eqref{randlin} yields the assertion of the 
proposition.
\hfill $\square$

\bigskip\noindent
By virtue of Proposition \ref{hardexp}, the hard constraints in \eqref{firstopt} define a polyhedral constraint set for the decision vector $x$. 
Recalling Lemma \ref{betconv}, \eqref{firstopt} would be a convex optimization problem provided that the probabilistic constraint defines a convex feasible region. As 
discussed in Section \ref{proconldr}, this can be guaranteed, however, only in certain special cases. Moreover, the range of applicability of Proposition \ref{hardexp} is potentially small:
\begin{cor}\label{hardexpcor}
Assume that all coefficients $\theta_{t, k}$ in \eqref{eqxtm} have all components different from zero. Then, if the hard constraints in \eqref{firstopt} represent simple box constraints, the only feasible linear decision rules are static ones.
\end{cor}
{\textbf{Proof.}}
For box constraints $y\in [y^{lo},y^{up}]$, we are dealing with the data specified in \eqref{boxexp}. Accordingly, the equation derived in Proposition \ref{hardexp} yields that
\[
A_{t, t}^{(3)}F_t\Theta_{1:t-1}+B_{t, t}^{(3)}\Theta_t=0\quad t=1,\ldots ,T.
\]
Recalling that, by the assumed wait-and-see structure for the hard constraints, we have $B_{t, t}^{(3)}=0$ for $t=1,\ldots ,T$ (see 
Section \ref{hardwait}), and taking into account that $A_{t, t}^{(3)}=(I,-I)^T$, we derive in particular 
the relations $F_t\Theta_{1:t-1}=0$ for $t=1,\ldots ,T$. Now, our assumption on coefficients 
$\theta_{t,k}$ 
ensures that the matrices $\Theta_{1:t-1}$ are surjective. As a consequence, $F_t=0$ for $t=1,\ldots ,T$, which means that the linear decision rules in \eqref{ldr} reduce to $y_{t}\left( \xi _{1:t-1}\right)=f_t$ for $t=1,\ldots ,T$. 
In other words, one is back to a static decision problem.\hfill $\square$

\bigskip\noindent
In order to avoid the restrictive consequences following from the last 
corollary, one may pass from Gaussian to truncated Gaussian distributions having a bounded support. This will be discussed in Section \ref{truncgaussopt}.
\subsection{Second optimization problem}
The second optimization problem to be discussed is \eqref{opt2}. We will focus our attention on the inner optimization problem 
\begin{equation}\label{inneropt}
\min\{h(y)\mid y\in M_2\cap\mathcal{K}\}.
\end{equation}
If this problem happens to have a unique solution, then its projection via $\Pi$ onto the hard constraints will be unique 
 and thus will be a solution of the overall problem too. Otherwise, the outer optimization problem in \eqref{opt2} just serves the purpose of selecting the best solution among projected solutions of the inner problem possibly realizing different values of the objective function $h$. We will not address the issue of possible non-uniqueness of \eqref{inneropt} here.

By \eqref{m2def}, and using once more the compact notation of Section \ref{lindec} along with the definition \eqref{ldr} of linear decision rules,  problem \eqref{inneropt} writes
\begin{equation}\label{secondopt}
\min\{\mathcal{J}(x)\mid\mathbb{P}(G^{(2)}_t(x)\varepsilon\leq g^{(2)}_t(x),\,\,G^{(3)}_t(x)\varepsilon\leq g^{(3)}_t(x)\quad (t=1,\ldots ,T))\geq p\}.
\end{equation}
This problem has the same objective as \eqref{firstopt} but the feasible set differs by the absence of hard constraints and the presence of an enlarged inequality system in the 
joint chance constraint. Once, a solution $x^*$ of \eqref{secondopt} has been determined, it is projected onto the hard constraints (either using an explicit formula if possible or by solving a quadratic optimization problem as described in Section \ref{projldr}) in order to yield a decision policy $\Pi (x^*)$ which is feasible for the original 
infinite-dimensional problem \eqref{finopt}.
\subsection{Third optimization problem} \label{thirdoptsection}
The third optimization problem we consider is \eqref{opt3} or its equivalent form \eqref{opt5}. Observe first, that \eqref{opt3} can be written 
\[
\min \{h(z)|z=\Pi (y),\,y\in M_{2}\cap \mathcal{K}\}.
\]
The inclusion in the constraint set of this optimization problem is the same as in \eqref{inneropt} and can thus be formulated as the probabilistic constraint in \eqref{secondopt} under our convention  $x:=(F_t,f_t)_{t=1,\ldots T}$ . 
Taking into account formula \eqref{quadropt} for the projection $z=\Pi (y)$, we arrive at the following description for problem \eqref{opt3}:
\begin{eqnarray}\label{bilev}
\min\{h(z)&\mid& z_t(\xi_{1:t-1})=\underset{u}{\rm argmin}\{\vartheta_t (x,u,\xi)\mid\gamma_t (u,\xi)\leq 0\quad \forall\xi,\,\,\forall t=1,\ldots ,T\},\\
&&\mathbb{P}(G^{(2)}_t(x)\varepsilon\leq  g^{(2)}_t(x),\,\,G^{(3)}_t(x)\varepsilon\leq g^{(3)}_t(x)\quad (t=1,\ldots ,T))\geq p\},\nonumber
\end{eqnarray}
where 
\begin{eqnarray*}
\vartheta_t (x,u,\xi)&:=&\|F_t\xi _{1:t-1}+f_t-u\|^2\\
\gamma_t (u,\xi)&:=&A_{t, t}^{(3)}u+\sum\limits_{\tau =1}^{t-1}B_{t, \tau }^{(3)}\xi _{\tau
}+\sum\limits_{\tau =1}^{t-1}A_{t, \tau }^{(3)}z_{\tau }(\xi_{1:\tau -1})-
b_{t}^{(3)}
\end{eqnarray*}
(recall that due to successive resolution of constraints in \eqref{quadropt} the terms 
$z_{\tau }(\xi_{1:\tau -1})$ are known in step $t$ for $\tau=1,\ldots ,t-1$). Formally, 
\eqref{bilev} represents a kind of bilevel problem in variables $(x,z)$, where the upper-level variable $x$ is subjected to a joint probabilistic constraint and the lower-level variable $z$ is subjected to a continuum of 
lower-level problems depending on $x$. As such, this optimization problem appears to be very hard to solve. On the other hand, for given $x$ satisfying the probabilistic constraint, the solutions $z_t$ of the parametric 
lower-level quadratic problem are piecewise linear in $\xi_{1:t-1}$ with an identifiable 
polyhedral decomposition of their domain. This would allow us to apply algorithms from multiparametric quadratic programming (see \cite{tondel}) in order to determine the $z_t$.

The problem simplifies significantly if the hard constraints are simple box constraints \eqref{boxhard} such that the explicit formula \eqref{boxform} can be applied. In this case, one may directly pass to the equivalent problem \eqref{opt5} which in our compact notation reads
\begin{eqnarray}\label{thirdopt}
&\mbox{minimize}&\nonumber\\
&\sum\limits_{t=1}^{T}\mathbb{E}\left\{ \langle h_{t},\delta_{t}(x,\xi ) \rangle + \left< \mathcal{P}_{t}, \left( \sum\limits_{\tau
=1}^{t}A_{t, \tau }^{(1)} \delta_{\tau }(x,\xi )+\sum\limits_{\tau
=1}^{t}B_{t, \tau }^{(1)}\xi _{\tau }-b_{t}^{(1)}\right) _{+} \right> \right\}\\
&\mbox{subject to}&\nonumber\\
&\mathbb{P}(G^{(2)}_t(x)\varepsilon\leq g^{(2)}_t(x),\,\,G^{(3)}_t(x)\varepsilon\leq g^{(3)}_t(x)\quad (t=1,\ldots ,T))\geq p,&\nonumber
\end{eqnarray}
where $x:=(F_t,f_t)_{t=1,\ldots T}$ and the components of $\delta_{t}(x,\xi )$ are defined as
\begin{equation}\label{projlindec}
(\delta_{t}(x,\xi ))_i:=\left(\max\{(\underline{y}_t)_i,\min\left\{(F_t\xi _{1:t-1}+f_t)_i,(\overline{y}_t)_i\right\}\}\right)_{i;\,t=1,\ldots ,T}.
\end{equation}
The first part of the expectation in the objective of this problem requires just to compute the expectations $\mathbb{E}(\delta_{t}(x,\xi ))_i$ which can be made fully explicit
thanks to Lemma \ref{lemmeanpiece} upon putting there (see \eqref{substi2})
\[
a:=(\underline{y}_t)_i;\,\,b:=(\overline{y}_t)_i;\,\,m:=(F_t\tilde{\mu}_{1:t-1}+f_t)_i;\,\,\sigma:=\sqrt{\left(F_t\Theta_{1:t-1}\Sigma\Theta_{1:t-1}^TF_t^T\right)_{ii}}.
\]
Consequently, in the absence of penalty terms in the objective, the whole problem reduces to a standard optimization problem subject to joint linear probabilistic constraints with multivariate Gaussian distribution. It may be difficult to obtain an analytic expression for the expectation of the penalty terms applied to projected linear decision rules. In this case, more elementary techniques like 
Sample Average Approximation may be used to approximate these expectations numerically.
\subsection{Fourth optimization problem}
The last optimization problem we consider is \eqref{opt4}. The difference with the previous op\-ti\-mization problems is that here decision variables are projections onto hard constraints from the very beginning. Similarly to the previous optimization problem, \eqref{opt4} can be written
\begin{equation}\label{fourthopt}
\min \{h(z)|z=\Pi (y),\,z\in M_1,\,y\in\mathcal{K}\}.
\end{equation}
Since the projection $z=\Pi (y)$ already ensures the hard constraint in the inclusion $z\in M_1$, it is sufficient to impose the probabilistic constraint in \eqref{m1def} on $z$. Following the idea and the notation of \eqref{bilev} in the previous optimization problem, one may reformulate \eqref{opt4} as
\begin{eqnarray}\label{bilev2}
\min\{h(z)&\mid& z_t(\xi_{1:t-1})=\underset{u}{\rm argmin}\{\vartheta_t (x,u,\xi)\mid\gamma_t (u,\xi)\leq 0\quad \forall\xi,\,\,\forall t=1,\ldots ,T\},\\
&&\mathbb{P}\left( \sum\limits_{\tau =1}^{t}A_{t, \tau }^{(2)}z_{\tau }\left(
\xi _{1:\tau -1}\right) +\sum\limits_{\tau =1}^{t}B_{t, \tau }^{(2)}\xi _{\tau
}\leq b_{t}^{(2)},\quad t=1,\ldots ,T\right) \geq p\}.\nonumber
\end{eqnarray}
Again, we are dealing with a bilevel problem in variables $(x,z)$, where the lower-level variable $z$ is subjected to a continuum of lower-level problems depending on the upper-level variable $x$. This time, however, the probabilistic constraint does not operate on the upper but rather on the 
lower-level variable. Moreover, it involves only the system of soft constraints (labeled by the upper index '(2)'). Evidently, in solving \eqref{bilev2} one is faced with the same difficulties as for problem \eqref{bilev}.

As before, there is motivation to investigate the special case of box constraints \eqref{boxhard}. Since in this case the projection $\Pi (y)$ can be made explicit via \eqref{boxform}, we may equivalently write \eqref{fourthopt}
\[
\min \{h(\Pi (y))\mid\Pi (y)\in M_1,\,y\in\mathcal{K}\}.
\]
This problem has the same objective as problem \eqref{opt5} and, hence, can be made explicit exactly the same way as described in the previous section for \eqref{thirdopt}. The difference now comes with the 
occurrence of {\it projected} linear decision rules \eqref{projlindec} as variables in the probabilistic constraint of \eqref{bilev2}. More precisely, we are led to the following optimization problem (where again $x:=(F_t,f_t)_{t=1,\ldots T}$):
\begin{eqnarray}\label{fourthopt1}
&\mbox{minimize}&\nonumber\\
&\sum\limits_{t=1}^{T}\mathbb{E}\left\{ \langle h_{t},\delta_t(x,\xi ) \rangle + \left< \mathcal{P}_{t}, \left( \sum\limits_{\tau
=1}^{t}A_{t, \tau }^{(1)} \delta_{\tau }(x,\xi )+\sum\limits_{\tau
=1}^{t}B_{t, \tau }^{(1)}\xi _{\tau }-b_{t}^{(1)}\right) _{+} \right> \right\}\\
&\mbox{subject to}&\nonumber\\
&\mathbb{P}\left( \sum\limits_{\tau =1}^{t}A_{t, \tau }^{(2)}\delta_{\tau}(x,\xi ) +\sum\limits_{\tau =1}^{t}B_{t, \tau }^{(2)}\xi _{\tau
}\leq b_{t}^{(2)},\quad t=1,\ldots ,T\right) \geq p.&\nonumber
\end{eqnarray}
The challenge now is to deal with the projected linear decision rules inside the probabilistic constraint and to reduce this issue to a tractable linear structure of type \eqref{pcfindim}. To this aim, with each index tuple
\[
(i_{1,1},\ldots ,i_{1,n_1},\ldots ,i_{T,1},\ldots i_{T,n_T})\in\{1,2,3\}^{\sum_{t=1}^Tn_t}
\]
we associate the following $x-$dependent partition of the space of events:
\begin{eqnarray*}
&S_{(i_{1,1},\ldots ,i_{1,n_1},\ldots ,i_{T,1},\ldots i_{T,n_T})}(x):=&\\ \\&\left\{\omega\in\Omega\mid\left\{\begin{array}{ll}
(F_t\xi _{1:t-1}(\omega)+f_t)_j\leq(\underline{y}_t)_j&\mbox{if }i_{t,j}=1\\
(\underline{y}_t)_j\leq (F_t\xi _{1:t-1}(\omega)+f_t)_j\leq (\overline{y}_t)_j&\mbox{if }i_{t,j}=2\\
(F_t\xi _{1:t-1}(\omega)+f_t)_j\geq(\overline{y}_t)_j&\mbox{if }i_{t,j}=3\\
\end{array}\right.\right\}&.
\end{eqnarray*}
Actually, this not a partition in the strict sense because the case distinction in its definition allows some overlap for nonstrict inequality signs. Due to $\xi$ having a density, however, this overlap is of measure zero. Therefore, we are allowed to reformulate the probability function in \eqref{fourthopt} as
\[
\sum\limits_{(i_{1,1},\ldots ,i_{1,n_1},\ldots ,i_{T,1},\ldots i_{T,n_T})\in\{1,2,3\}^{\sum_{t=1}^Tn_t}}\mathbb{P}\left(
\begin{array}{l}
\xi\in S_{(i_{1,1},\ldots ,i_{1,n_1},\ldots ,i_{T,1},\ldots i_{T,n_T})}(x),\\
\sum\limits_{\tau =1}^{t}\sum\limits_{j=1}^{n_\tau}(\delta_{\tau}(x,\xi ))_j(A_{t, \tau }^{(2)})_j +\sum\limits_{\tau =1}^{t}B_{t, \tau }^{(2)}\xi _{\tau
}\leq b_{t}^{(2)}\\
(t=1,\ldots ,T)
\end{array}
\right),
\]
where $(A_{t, \tau }^{(2)})_j$ refers to column $j$ of the matrix $A_{t, \tau }^{(2)}$. Observing that, by definition,
\[
(\delta_{\tau}(x,\xi ))_j=\left\{\begin{array}{ll}
(\underline{y}_\tau)_j&\mbox{if }i_{\tau,j}=1\\
(F_\tau\xi _{1:\tau -1}+f_\tau)_j&\mbox{if }i_{\tau,j}=2\\
(\overline{y}_\tau)_j&\mbox{if }i_{\tau,j}=3
\end{array}\right.,
\]
we realize that each event over which the probability is taken above, is described by a system of random inequalities which is linear in the random vector $\xi$. Consequently, the probability of each such event above can be described by
\[
\mathbb{P}\left(\tilde{G}_t^{(i_{1,1},\ldots ,i_{1,n_1},\ldots ,i_{T,1},\ldots i_{T,n_T})}(x)\xi\leq\tilde{g}_t^{(i_{1,1},\ldots ,i_{1,n_1},\ldots ,i_{T,1},\ldots i_{T,n_T})}(x)
\quad (t=1,\ldots ,T)\right).
\]
With $\xi$ being an affine linear mapping of $\varepsilon$ according to \eqref{substi2}, we may finally write the probabilistic 
constraint in \eqref{fourthopt} as
\[
\sum\limits_{\begin{array}{c}(i_{1,1},\ldots ,i_{1,n_1},\ldots ,i_{T,1},\ldots i_{T,n_T})\\\in\{1,2,3\}^{\sum_{t=1}^Tn_t}\end{array}}\mathbb{P}\left(
\begin{array}{l}\tilde {G}_t^{(i_{1,1},\ldots ,i_{1,n_1},\ldots ,i_{T,1},\ldots i_{T,n_T})}(x)\xi\leq\\ {\tilde g}_t^{(i_{1,1},\ldots ,i_{1,n_1},\ldots ,i_{T,1},\ldots i_{T,n_T})}(x)\quad (t=1,\ldots ,T)
\end{array}\right)\geq p,
\]
which now involves similar terms as \eqref{pcfindim}. 

Clearly this approach for dealing with the probabilistic constraint in \eqref{fourthopt}
quickly becomes prohibitive due to the number $3^{\sum_{t=1}^Tn_t}$ of terms in the sum above. Even if every decision policy is one-dimensional ($n_t=1$ for all $t$), this yields $3^T$ summands and limits the applicability of the approach to say $T=6,7$ stages. An 
alternative option would consist in the application of spherical-radial decomposition as mentioned in Section \ref{proconldr} which is not restricted to linear probabilistic constraints and would not suffer from the complexity issue.
\subsection{Optimization problem under truncated Gaussian distribution}\label{truncgaussopt}
After introducing our original optimization problem \eqref{finopt}, we have passed 
immediately to hard constraints of wait-and-see type in Section \ref{hardwait} because otherwise the hard constraints would not have any good chance of ever being satisfied under distributions with unbounded support, e.g., Gaussian. This change became apparent by 
requiring $B_{t, t}^{(3)}=0$ in \eqref{finopt}, leading to the hard constraints of \eqref{m1def}. When discussing our first 
optimization problem \eqref{firstopt}, we noticed that even for hard constraints of wait-and-see type, the unboundedness of 
the support of the random vector generates a strong restriction on the feasible decisions (see Corollary \ref{hardexpcor}). In this section we come back to the first optimization problem but with a Gaussian random vector truncated to a bounded region. This approach will allow us not only to circumvent the mentioned restriction of problem \eqref{firstopt} but even to admit the original hard constraints in \eqref{finopt} with possibly $B_{t, t}^{(3)}\neq 0$. 
\begin{df}
We say that a random vector $\varepsilon$ follows a normal distribution with parameters $\mu,\Sigma$ which is truncated to a Borel measurable set $S$ and then write $\varepsilon\sim\mathcal{TN} (\mu,\Sigma,S)$ if there exists a Gaussian random vector $\tilde{\varepsilon}\sim\mathcal{N} (\mu,\Sigma)$ such that
\[
\mathbb{P}(\varepsilon\in B)=\frac{\mathbb{P}(\tilde{\varepsilon}\in S\cap B)}{\mathbb{P}(\tilde{\varepsilon}\in S)}\quad\mbox{for all Borel sets }B.
\]
\end{df}
In the following we shall assume in contrast with the previous sections that the noises 
$\varepsilon_t$ in the probabilistic model \eqref{modeletimeseries} are independent and
distributed according 
to $\varepsilon\sim\mathcal{TN} (0,\Sigma,S)$, where $\Sigma$ is the block-diagonal matrix introduced in 
Section \ref{probmod}.

\if{
{\color{red}\noindent Je ne sais pas si on a besoin du lemme suivant. Est-ce qu'on a besoin dans la th\'eorie des s\'eries temporelles que les innovations al\'eatoires soient ind\'ependantes? Dans le papier original - sous (18) - tu mentionnes que les composantes de $\varepsilon$ sont ind\'ependantes. Je pense que cela n'est pas automatique parce que tout ce qu'on a jusque-l\`a, c'est l'ind\'ependance des composantes du Gaussien $\tilde{\varepsilon}$ gr\^ace \`a la block-diagonalit\'e de $\Sigma$. Par contre, il n'est pas n\'ec\'essaire de mettre l'ind\'ependance des composantes de $\varepsilon$ comme une hypoth\`ese parce qu'elle suit de celle des composantes de $\tilde{\varepsilon}$.}
\begin{lem}
The components of the random vector $\varepsilon=(\varepsilon_1,\ldots ,\varepsilon_T)$ are independent.
\end{lem}
{\textbf{Proof.}}
Let $\tilde{\varepsilon}\sim\mathcal{N} (\mu,\Sigma)$ and write the box $S=S_1\times\cdots\times S_T$ as a product of boxes. Fix any couple $r,s\in\{1,\ldots ,T\}$ with $r\neq s$. Then, for any couple of Borel sets $B_r,B_s$ it holds 
\begin{eqnarray*}
\mathbb{P}(\varepsilon_r\in B_r,\,\varepsilon_s\in B_s)&=&
\frac{\mathbb{P}(\tilde{\varepsilon}_r\in B_r\cap S_r,\,\tilde{\varepsilon}_s\in B_s\cap S_s,\,\tilde{\varepsilon}_t\in S_t\,\forall t\in\{1,\ldots ,T\}\setminus\{r,s\}}{\mathbb{P}
(\tilde{\varepsilon}\in S)}\\
&=&\frac{\mathbb{P}(\tilde{\varepsilon}_r\in B_r\cap S_r)\cdot\mathbb{P}(\tilde{\varepsilon}_s\in B_s\cap S_s)\cdot\Pi_{t\in\{1,\ldots ,T\}\setminus\{r,s\}}\mathbb{P}(\tilde{\varepsilon}_t\in S_t)}{\mathbb{P}(\tilde{\varepsilon}_r\in S_r)\cdot\mathbb{P}(\tilde{\varepsilon}_s\in S_s)\cdot\Pi_{t\in\{1,\ldots ,T\}\setminus\{r,s\}}\mathbb{P}(\tilde{\varepsilon}_t\in S_t)}\\
&=&\mathbb{P}(\varepsilon_r\in B_r)\cdot\mathbb{P}(\varepsilon_s\in B_s).
\end{eqnarray*}
Here, we exploited that the components of the random vector $\tilde{\varepsilon}=(\tilde{\varepsilon}_1,\ldots ,\tilde{\varepsilon}_T)$ are independent by block diagonality of $\Sigma$.

\hfill
$\square$
}\fi
We are now going to check the impact of truncating the Gaussian distribution on the structure of optimization problem \eqref{firstopt}. 

\if{
{\color{red}\noindent Dans le papier original au milieu de la page 12 tu \'ecris que l'objectif du probl\`eme est le m\^eme que celui dans (7) du papier original. Mais maintenant $\varepsilon$ a chang\'e de distribution, il n'est m\^eme pas clair que l'esp\'erance de $\varepsilon$ soit z\'ero (sauf que si $S$ est un ensemble centr\'e). \`A mon avis, cela change le tableau dans le calcul de l'objectif. Qu'est-ce que tu penses?} 
}\fi
The terms 
$\mathbb{E}[\langle h_{t},y_{t}\left( \xi_{1:t-1}\right) \rangle ]$ in the objective function
can be computed analytically since closed-form expressions are available for the expectation
of truncated normal one-dimensional random variables.

\if{
\begin{lem} Let $X$ be a one dimensional random variable
distributed according to $X \sim\mathcal{TN} (m,\sigma^2,[a,b])$ with $a<b$,
let $f(x)=\frac{1}{\sqrt{2 \pi}} e^{-\frac{x^2}{2}}$ be the density of the
standard Gaussian distribution, and let $\Phi$ be
the one-dimensional standard normal distribution function.
Then
$$
\mathbb{E}[X]=m + \left(\frac{f\left(\frac{a-m}{\sigma}\right)-f\left(\frac{b-m}{\sigma}\right)}{\Phi\left(\frac{b-m}{\sigma} \right)- \Phi\left(\frac{a-m}{\sigma}\right)}\right) \sigma.
$$
\end{lem}
}\fi
Similarly to problem \eqref{thirdopt}, the expectation of the penalty terms can be approximated using Sample Average Approximation.

If $S:=[\underline{S},\overline{S}]$ is a box, as far as the probabilistic constraint in \eqref{firstopt} is concerned, the underlying probability function can be written 
\begin{eqnarray*}
\mathbb{P}(G_t^{(2)}(x) \varepsilon \leq g_t^{(2)}(x)\,\,(t=1,\ldots ,T))&=&\\\frac{ \mathbb{P}(\{G_t^{(2)}(x) {\tilde \varepsilon} \leq g_t^{(2)}(x)\,\,(t=1,\ldots ,T)\}  \cap  \{ {\tilde \varepsilon} \in S\})}{\mathbb{P}({\tilde \varepsilon} \in S)}&=&\frac{\mathbb{P}({\tilde G}(x) {\tilde \varepsilon} \leq {\tilde g}(x))} {\mathbb{P}({\tilde \varepsilon} \in S)},
\end{eqnarray*}
where, with $I$ referring to the identity matrix of appropriate size, 
\[
\tilde{G}(x):=
\left(
\begin{array}{c}
G_1^{(2)}(x)\\
\vdots\\ 
G_T^{(2)}(x)\\
I\\
-I
\end{array}
\right),\quad\tilde{g}(x):=
\left(
\begin{array}{c}
g_1^{(2)}(x)\\
\vdots\\
g_T^{(2)}(x)\\
\overline{S}\\
-\underline{S}
\end{array}
\right).
\]
Consequently, the probabilistic constraint in \eqref{firstopt} turns into
\begin{equation}\label{reformjcclognorm}
\mathbb{P}({\tilde G}(x) {\tilde \varepsilon} \leq {\tilde g}(x))\geq\tilde{p},\quad
\mbox{where }\tilde{p}:=p\cdot\mathbb{P}({\tilde \varepsilon} \in S).
\end{equation}
Due to $\tilde{\varepsilon}$ being a Gaussian random vector, this probabilistic constraint is exactly of the same nature as the original one in \eqref{firstopt} which was discussed in Section \ref{proconldr}.

Addressing finally the almost sure constraints in \eqref{firstopt}, they can be equivalently formulated as
\begin{equation} \label{robustconstraints}
\max\limits_{\varepsilon \in S} \left(G_t^{(3)}(x)\right)^j \varepsilon\leq g_{t,j}^{(3)}(x),\quad\forall t,\,\,\forall j,
\end{equation}
where $\left(G_t^{(3)}(x)\right)^j$ refers to the $j$th line of $G_t^{(3)}(x)$. 

We consider two cases for $S$: a box and an ellipsoid.
 If $S:=[\underline{S},\overline{S}]$ is a box, the maximum in the 
left-hand side of \eqref{robustconstraints} can be computed analytically using the following lemma:
\begin{lem}[\cite{robprodman}, Lemma 2]
For any $x$ we have that
\[
\max_{y \in S}\;x^T y=\frac{1}{2} \left(x^T(\underline{S}+\overline{S})+|x|^T (\overline{S}-\underline{S}) \right).
\]
\end{lem}
As a result, if $S$ is a box, since $\left(G_t^{(3)}(x)\right)^j$ and $h_{t,j}^{(3)}(x)$ are affine functions of $x$,
the almost sure constraints in \eqref{firstopt} 
can be reformulated as explicit convex constraints in $x$. 
\if{
Setting $S_c=(\underline{S} +  \overline{S})/2$ and $S_L = (\overline{S} - \underline{S})/2 \geq 0$,
we can define the following linear relaxation of constraint \eqref{robustconstraints}:
$$
\left\{
\begin{array}{l}
S_c^T  \left(G_t^{(3)}(x)\right)^j + S_L^T  z_{t, j} \leq  h_{t,j}^{(3)}(x),\\
z_{t,j}^T  \geq \left( G_t^{(3)}(x)\right)^j, \,z_{t,j}^T \geq -\left( G_t^{(3)}(x)\right)^j,  
\end{array}
\right.
$$
where $z_{t, j}$ is a vector in $\mathbb{R}^T$.
}\fi

Now taking for $S$ the ellipsoid
$$
S=\{x \in \mathbb{R}^T : (x - \mu)^T \Sigma^{-1} (x - \mu) \leq \kappa^2   \},
$$
if vector $w_{t,j}(x)$ is the transpose of $\left(G_t^{(3)}(x)\right)^j$ then
constraint \eqref{robustconstraints} can be reformulated as the explicit conic quadratic (convex) constraint
$$
\mu^T w_{t, j}(x) + \kappa \sqrt{ w_{t, j}(x)^T  \Sigma   w_{t, j}(x)    } \leq h_{t,j}^{(3)}(x).
$$
We end up again with a convex optimization problem. 

Finally, observe that the term $\mathbb{P}({\tilde \varepsilon} \in S)$ in \eqref{reformjcclognorm}
can be computed numerically when $S$ is a box (using Genz' code as described in \cite{genz} for instance).

\section*{Acknowledgments} The first author's research was supported by an FGV grant, CNPq grant 307287/2013-0, 
FAPERJ grants E-26/110.313/2014, and E-26/201.599/2014.
The second author gratefully acknowledges support by the {\it  FMJH Program Gaspard Monge in optimization and operations research} including support to this program by 
EDF as well as support by the {\it Deutsche Forschungsgemeinschaft} within Projekt B04 in CRC TRR 154.

\bibliographystyle{plain}
\bibliography{DJCC}

\section*{Appendix}

\section{Algorithm for computing coefficients $c, \gamma, \delta$, and $\theta$  of decomposition \eqref{eqxtm}}
\begin{lem}\label{lemcalcijh} Let $\xi_t$ satisfy \eqref{modeletimeseries} and for any positive integers $t, j$, let $I_{t, j}, J_{t, j}$,
and $H_{t, j}$ be the sets given by
\begin{eqnarray*}
I_{t, j}(m) &= & \{k \in \mathbb{N} \;:\; 1 \leq k \leq \min(p_{t+1}(m),t+1-j)\},\\
J_{t, j}(m) &= & \{k \in \mathbb{N} \;:\; 1 \leq k \leq \min(t,p_{t+1}(m)),\;j \leq r_{t+1-k}(m)\},\\
H_{t, j}(m) &= & \{k \in \mathbb{N} \;:\; 1 \leq k \leq \min(t,p_{t+1}(m)),\;j \leq s_{t+1-k}(m)\}.
\end{eqnarray*}
We also define 
\begin{eqnarray*}
X_t(m) & = & \max\left(r_{t+1-k}(m), k=1,\ldots,\min(t,p_{t+1}(m))\right),\\ 
Y_t(m) & = & \max\left(s_{t+1-k}(m), k=1,\ldots,\min(t,p_{t+1}(m))\right).
\end{eqnarray*}
The coefficients $c, \gamma, \delta,$ and $\theta$ in the decompositions
of $\xi_1(m), \xi_2(m)$, $\ldots,\xi_T(m), m=1,\ldots,M$, of the form \eqref{eqxtm} are computed iteratively
as follows:\\
\par {\underline{Initialization:}} For $m=1,\ldots,M$, set $c_1(m)=\frac{\mu_1(m)}{\alpha_{1, 0}(m)}$, 
$r_1(m)=p_1(m)$, $\gamma_{1, k}(m)=-\frac{\alpha_{1, k}(m)}{\alpha_{1, 0}(m)}, k=1,\ldots,p_1(m)$,
$s_1(m)=q_1(m)$, $\delta_{1, k}(m)=\frac{\beta_{1, k}(m)}{\alpha_{1, 0}(m)}, k=1,\ldots,q_1(m)$, and 
$\theta_{1, 1}(m)=\frac{\beta_{1, 0}(m)}{\alpha_{1, 0}(m)}$.\\
\par {\underline{Loop:}} For $m=1,\ldots,M$ and for $t=1,\ldots,T-1$,
$$
c_{t+1}(m)=\frac{\mu_{t+1}(m)}{\alpha_{t+1, 0}(m)}-\displaystyle{\sum_{k=1}^{\min(t, p_{t+1}(m))}} \frac{\alpha_{t+1, k}(m)} {\alpha_{t+1, 0}(m)} c_{t+1-k}(m).
$$
$$
\theta_{t+1, j}(m)=\left\{
\begin{array}{l}
\frac{\beta_{t+1, 0}(m)}{\alpha_{t+1, 0}(m)}  \;\emph{ for }\; j=t+1,\\
\frac{\beta_{t+1, t+1-j}(m)}{\alpha_{t+1, 0}(m)}-\displaystyle{\sum_{k \in I_{t, j}(m)}} \frac{\alpha_{t+1, k}(m)}{\alpha_{t+1, 0}(m)} \theta_{t+1-k, j}(m)  \;\emph{ for }\;t+1-\min(t, q_{t+1}(m)) \leq j \leq t,\\
-\displaystyle{\sum_{k \in I_{t, j}(m)}} \frac{\alpha_{t+1, k}(m)}{\alpha_{t+1, 0}(m)} \theta_{t+1-k, j}(m)  \;\emph{ for }\;1 \leq j \leq t-\min(t, q_{t+1}(m)).
\end{array}
\right.
$$
Coefficient $\gamma_{t+1, j}(m)$ is given by
$$
\left\{
\begin{array}{l}
-\frac{\alpha_{t+1, j+t}(m)}{\alpha_{t+1, 0}(m)}-\displaystyle{\sum_{k \in J_{t, j}(m)}} \; \frac{\alpha_{t+1, k}(m)}{\alpha_{t+1, 0}(m)} \gamma_{t+1-k, j}(m)  \;\emph{ for }\;1 \leq j \leq \min(p_{t+1}(m)-t, X_t(m)),\\
-\frac{\alpha_{t+1, t+j}(m)}{\alpha_{t+1, 0}(m)}   \;\emph{ for }\;1+\min(p_{t+1}(m)-t,X_t(m)) \leq j \leq p_{t+1}(m)-t,\\
-\displaystyle{\sum_{k \in J_{t, j}(m)}} \frac{\alpha_{t+1, k}(m)}{\alpha_{t+1, 0}(m)} \gamma_{t+1-k, j}(m)  \;\emph{ for }\;\max\left(1,1+\min(p_{t+1}(m)-t, X_t(m))\right) \leq j \leq X_t(m).
\end{array}
\right.
$$
Coefficient $\delta_{t+1, j}(m)$ is given by
$$
\left\{
\begin{array}{l}
\frac{\beta_{t+1, j+t}(m)}{\alpha_{t+1, 0}(m)}-\displaystyle{\sum_{k \in H_{t, j}(m)}} \; \frac{\alpha_{t+1, k}(m)}{\alpha_{t+1, 0}(m)} \delta_{t+1-k, j}(m)  \;\emph{ for }\;1 \leq j \leq \min(q_{t+1}(m)-t, Y_t(m)),\\
\frac{\beta_{t+1, t+j}(m)}{\alpha_{t+1, 0}(m)}   \;\emph{ for }\;1+\min(q_{t+1}(m)-t, Y_t(m)) \leq j \leq q_{t+1}(m)-t,\\
-\displaystyle{\sum_{k \in H_{t, j}(m)}} \frac{\alpha_{t+1, k}(m)}{\alpha_{t+1, 0}(m)} \delta_{t+1-k, j}(m)  \;\emph{ for }\;\max\left(1,1+\min(q_{t+1}(m)-t, Y_t(m))\right) \leq j \leq Y_t(m).
\end{array}
\right.
$$
Finally,
$$
r_{t+1}(m)=\max\left(p_{t+1}(m)-t, X_{t}(m)\right)  \mbox{ and } s_{t+1}(m)=\max\left(q_{t+1}(m)-t, Y_{t}(m)\right).
$$
\end{lem}
{\textbf{Proof.}}
We fix a component $m$ and to alleviate notation, we drop $(m)$ in the proof. The initialization is immediate,
writing \eqref{modeletimeseries} for $t=1$. Now assume that for some $t<T$, the decompositions 
of $\xi_1, \ldots, \xi_t$ of the form \eqref{eqxtm} are available. To obtain the decomposition of
$\xi_{t+1}$, we use \eqref{modeletimeseries} to obtain
\begin{eqnarray*}
\xi_{t+1} & = &  \frac{\mu_{t+1}}{\alpha_{t+1, 0}}-\displaystyle{\sum_{k=1}^{\min(t, p_{t+1})}} \frac{\alpha_{t+1, k}}{\alpha_{t+1, 0}} \xi_{t+1-k}+ \displaystyle{\sum_{k=0}^{\min(t, q_{t+1})}} \frac{\beta_{t+1, k}}{\alpha_{t+1, 0}} \varepsilon_{t+1-k}\\
& &  -\displaystyle{\sum_{k=1+\min(t, p_{t+1})}^{ p_{t+1}}} \frac{\alpha_{t+1, k}}{\alpha_{t+1, 0}} \xi_{t+1-k}+ 
\displaystyle{\sum_{k=1+\min(t, q_{t+1})}^{ q_{t+1}}} \frac{\beta_{t+1, k}}{\alpha_{t+1, 0}} \varepsilon_{t+1-k}.
\end{eqnarray*}
In the first sum, since for all index $k \in \{1, 2, \ldots, \min(t, p_{t+1})\}$ we have
$1 \leq t+1-k \leq t$, we know for $\xi_{t+1-k}$ 
a decomposition of the form \eqref{eqxtm} with known coefficients  $c, \gamma, \delta$, and $\theta$.
Using these expressions of $\xi_{t+1-k}$, this first sum can be written
$$
-\displaystyle{\sum_{k=1}^{\min(t, p_{t+1})}} \frac{\alpha_{t+1, k}}{\alpha_{t+1, 0}} \left(c_{t+1-k}  + \displaystyle{\sum_{j=1}^{r_{t+1-k}}} \gamma_{t+1-k, j} \xi_{1-j}+\displaystyle{\sum_{j=1}^{s_{t+1-k}}} \;\delta_{t+1-k, j} \varepsilon_{1-j} +\displaystyle{\sum_{j=1}^{t+1-k}} \;\theta_{t+1-k, j} \varepsilon_{j} \right).
$$
Gathering the terms that depend neither on noise $\varepsilon$ nor on $\xi$, we obtain the expression of
$c_{t+1}$.

The portion depending on $\varepsilon_1, \ldots, \varepsilon_{t+1}$ can be written
$$
\displaystyle{\sum_{j=t+1-\min(t, q_{t+1})}^{t+1}} \frac{\beta_{t+1, t+1-j}}{\alpha_{t+1, 0}} \varepsilon_j -\displaystyle{\sum_{j=1}^t} \left( \displaystyle{\sum_{k \in I_{t, j}}} \frac{\alpha_{t+1, k}}{\alpha_{t+1, 0}} \theta_{t+1-k, j}\right) \varepsilon_{j}.
$$
We then consider the decomposition of $\xi_{t+1}$ obtained replacing
$t$ by $t+1$ in \eqref{eqxtm}. Identifying the portion of this decomposition depending on $\varepsilon_1, \ldots, \varepsilon_{t+1}$
with the expression above, we obtain the expressions of the
coefficients $\theta_{t+1, j},\;j=1,\ldots, t+1$.

The portion that depends on $\xi_0, \xi_{-1},\ldots,$ can be written
$$
-\displaystyle{\sum_{j=1+\min(0, p_{t+1}-t)}^{p_{t+1}-t}} \frac{\alpha_{t+1, t+j}}{\alpha_{t+1, 0}} \xi_{1-j}-\displaystyle{\sum_{j=1}^{X_t}} \left( \displaystyle{\sum_{k \in J_{t, j}}} \frac{\alpha_{t+1, k}}{\alpha_{t+1, 0}} \gamma_{t+1-k, j}\right) \xi_{1-j}.
$$
From that expression, we obtain the desired value of $r_{t+1}$ as well as the announced formulas
for coefficients $\gamma_{t+1 j}, j=1,\ldots, r_{t+1}$.

Finally, the portion  depending on $\varepsilon_0, \varepsilon_{-1}, \ldots$, can be written
$$
\displaystyle{\sum_{j=1+\min(0, q_{t+1}-t)}^{q_{t+1}-t}} \frac{\beta_{t+1, t+j}}{\alpha_{t+1, 0}} \varepsilon_{1-j}-\displaystyle{\sum_{j=1}^{Y_t}} \left( \displaystyle{\sum_{k \in H_{t, j}}} \frac{\alpha_{t+1, k}}{\alpha_{t+1, 0}} \delta_{t+1-k, j}\right) \varepsilon_{1-j}.
$$
From that expression, we obtain the desired value of $s_{t+1}$ as well as the announced formulas
for coefficients $\delta_{t+1, j}, j=1,\ldots, s_{t+1}$.\hfill $\square$

\end{document}